\documentclass[a4paper,12pt]{article}
\usepackage{amsmath}
\usepackage{amssymb}
\newcommand{\reals}{\mathbb {R}}
\newcommand{\integers}{\mathbb{Z}}

\newcommand{\complexes}{\mathbb{C}}

\newcommand{\So}[1]{{\mathrm {SO}}(#1)}

\renewcommand{\O}[1]{{\mathrm O}(#1)}
\renewcommand{\o}[1]{{\mathfrak o}(#1)}
\newcommand{\U}[1]{{\mathrm U}(#1)}
\renewcommand{\u}[1]{{\mathfrak u}(#1)}
\newcommand{\Su}[1]{{\mathrm {SU}}(#1)}

\newcommand{\SP}[1]{{\mathrm {Sp}}(#1)}
\newcommand{\sP}[1]{{\mathfrak {sp}}(#1)}

\newcommand{\Hom}[2]{{\mathrm {Hom}}(#1,#2)}
\newcommand{\End}[1]{{\mathrm {End}}(#1)}

\newcommand{\mapping}{\rightarrow}

\newcommand{\tensor}{\otimes}

\renewcommand{\mod}[1]{{\left | {#1} \right |}}

\newcommand{\iso}{\cong}

\newcommand{\grad}{\nabla}

\newcommand{\restrict}[1]{\bigg\vert_{#1}}
\newcommand{\inv}[1]{{#1}^{-1}}
\newcommand{\tang}[2]{{\mathrm{T}}_{#2}{#1}}
\newcommand{\Tang}[1]{{\mathrm{T}}{#1}}

\newcommand{\sphere}[1]{{\mathrm{S}}^{#1}}

\newcommand{\function}[3]{{#1}{\colon}{#2}{\mapping}{#3}}

\newcommand{\Symm}[2]{{\odot^{#2}}({#1})}

\newcommand{\Form}[2]{{\Omega^{#2}}(#1)}

\newcommand{\id}{\,\mathrm{l}\!\!\!\mathrm{1}}

\newcommand{\goth}[1]{\mathfrak #1}

 \newenvironment{prf}{{\bf Proof}\relax\\}{\hfill$\blacksquare$\par\bigskip}
 
 \newtheorem{defn1}{Definition}[section]
 \newtheorem{lem1}[defn1]{Lemma}
 \newtheorem{prop1}[defn1]{Proposition}
 \newtheorem{thm1}[defn1]{Theorem}
 \newtheorem{cor1}[defn1]{Corollary}

\usepackage{eucal}
\usepackage{epsfig}

\newcommand{\quat}{\mathbb{H}}

\newcommand{\sect}[2]{\Form{#1;#2}{0}}
\newcommand{\tE}{\hat{E}}

\newcommand{\Rf}{\reals^4}

\newcommand{\Diff}[1]{\frac{\partial}{\partial #1}}

\newcommand{\con}{\grad}

\newcommand{\I}{{\widehat{\textbf{i}}}}
\newcommand{\J}{{\widehat{\textbf{j}}}}
\newcommand{\K}{{\widehat{\textbf{k}}}}

\newcommand{\eps}{\varepsilon}

\renewcommand{\P}{\mathcal{P}}

\newcommand{\R}{\mathcal{R}}

\renewcommand{\mapping}{\longrightarrow}

\newcommand{\A}[1]{{\bf{A}^*}\left({#1}\right)}
\newcommand{\Abad}[1]{{\bf{A}}\left({#1}\right)}
\newcommand{\As}[1]{{\bf{A}}\left({#1}\right)}

\newcommand{\BBs}[1]{{M}_{#1}}
\newcommand{\BBt}[1]{\tilde{M}_{#1}}
\newcommand{\hcon}{\widehat{\grad}}
\newcommand{\hE}{\widehat{E}}
\newcommand{\OSP}[1]{\O{#1}\times\SP{1}}

\newcommand{\At}{\hbox{At}}
\renewcommand{\d}{\mathrm{d}}
\newcommand{\hd}{\widehat{\d}}
\newcommand{\F}{\mathcal{F}}

\newcommand{\hook}{
\begin{picture}(10,10)(-1,0)
\put(0,0){\line(1,0){7}}
\put(7,0){\line(0,1){7}}
\end{picture}}

\newcommand{\mc}[1]{{\mathfrak{M}_\complexes^{#1}}}
\newcommand{\mr}[1]{{\mathfrak{M}_\reals^{#1}}}

\newcommand{\Endoff}[1]{}

\newcommand{\wdot}{{\,.\!\!\!\wedge\,}}

\newcommand{\OG}[2]{\Omega_G^{#2}\left({#1}\right)}
\newcommand{\dg}{{\d_\goth{g}}}

\renewcommand{\exp}[1]{e^{{#1}}}
\newcommand{\g}{\goth{g}}

\newcommand{\conng}{\grad_\g}
\newcommand{\curv}{F_\g}

\title{The ADHM Construction and Anselmi's Topological Anomalies}
\author{Jonathan Munn}
\date{\today}
\begin{document}
%\cite{MU}
\maketitle
\begin{abstract}
We examine the anomalies arising in instanton calculus as detailed by Damiano
Anslemi in 1994. Whereas Anselmi uses BRST theory, we use the ADHM construction
to arrive at the same conclusions from a differential-geometric way. We observe
that Anselmi's TQFT is similar to Donaldson Theory applied to charge 1 instantons
on the 4-sphere,although the latter is really only used for 4-manifolds with
$b_2^+>0$. \par
We demonstrate why the anomalies occur in the case of charge 1 instantons
and move to show that similar anomalies cannot occur for instantons of higher
charge. To do this, we develop an equivariant integration theory for hyperKahler
manifolds and apply it to the hyperKahler geometry involved in the ADHM
construction.
\end{abstract}
\newpage
\tableofcontents
\section{Introduction}
In \cite{An}, Damiano Anselmi used BRST theory to discover certain topological
anomalies arising in one of the more simple examples of TFT, namely the
case of instantons, i.e bundles with ASD connections over $\sphere{4}$ (or
$\Rf$ with certain certain restrictions on the choice of gauge). The moduli
space of charge 1 instantons $\mathcal{M}_1$ (i.e bundles with $c_2=1$) of
course is equvialent to 5 dimensional hyperbolic space which is contractible,
hence all bundles over this space are trivial. As a result, all Donaldson polynomials
formed from the cohomology of $\mathcal{M}_1$ actually vanish in this case.\par
  Anselmi's method entailed building BRST representatives
of the Donaldson $\mu$ cohomology classes of $\mathcal{M}_1$ from data on the moduli space of charge 1, and integrating these over certain submanifolds
within $\mathcal{M}_1$. This process effectively calculates the Donaldson
polynomial associated with these submanifolds. However, Anselmi's approach produced a non-trivial linking theory in certain dimensions when integrals were taken over certain ``cycles''thus essentially contradicting the fact that there are no non-trivial Donaldson polynomials for $\sphere{4}$.\par
We use the ADHM construction as a basis for understanding the problem. This has the advantage that we can bypass BRST theory and obtain the same results
in which Anselmi's Anomalies manifest themselves as singularities of  Chern-Weil
representatives of characteristic classes of the moduli space and show that in essence the anomalies arise due to the noncommutative procedure of removing unstable points and hyperK\"ahler reduction. 
\par
We use an analogue of the Jeffrey-Kirwan localisation theorem for hyperK\"ahler
manifolds which relates the integration of equivariant cohomology classes
over a hyperK\"ahler manifold possessing a tri-Hamiltonian group action with the integration of cohomology classes on the quotient. We use this formula
to show that there are no further anomalies for higher charge
instantons essentially because of dimensional incompatibility.\par
%\input Section1.tex
%%\section{Preliminaries}
\subsection{Review of the ADHM construction}\label{Ansprelim}
From the work of Penrose and the famous paper \cite{ADHM}, it was shown that
vector bundles with ASD connections could be constructed over $\sphere{4}$ using the conformal relationship betwen $\sphere{4}$ and $\Rf$ and little more than quaternionic linear algebra. \par
Recall that we may form ASD $\Su{2}$-connections on $\Rf$ ($\sphere{4}$)
by choosing
\begin{eqnarray*}
(T,P)&\in
%& i\u{k}\tensor\Skew{\Cotang{\Rf}}{1}\oplus\Hom{\complexes^2}{\complexes^k\tensor\Splus}\\
%\end{eqnarray*}
%such that
%\[
%(T\wedge T)^++(PP^*)_0=0.
%\]
%Since, as real vector spaces
%\begin{eqnarray*}
%&i\u{k}\tensor\Skew{\Cotang{\Rf}}{1}\oplus\Hom{\complexes^2}{\complexes^k\tensor\Splus}\\
%&\iso
&\left(i\u{k}\tensor\quat\right)\oplus\left(\complexes^k\tensor_\reals\quat\right)\\
&\iso&\left(i\u{k}\oplus\complexes^k\right)\tensor_\reals\quat\\
&=:&\mc{k},
\end{eqnarray*}
%we can rewrite the ADHM condition as
such that
\[
\Im\left(T^*T+PP^*\right)=0
\]
where $\Im$ denotes the quaternionic imaginary part which is well defined as $\mc{k}$ is very much a ``quaternification'' of a real vector space.
We also require, for each $x\in\Rf$, that the map
\[
\R_x=\function{((T-x\id)^*,P)}{\complexes^k\tensor\quat\oplus\complexes^k}
{\complexes^k\tensor\quat}
\]
be surjective.
Then we define a bundle $E=\ker\R$ and a connection given by $v^*\hbox{d}v$,
where
\[
v_x=\left[
\begin{array}{c}
\inv{{(x\id-T)^*}}P \\
\id
\end{array}
\right]
\inv{\sigma_x}
\]
and
\[
{\sigma_x}^2=\id+P^*\inv{\left((T-x\id)^*(T-x)\right)}P.
\]
The quaternionic column vector $v$ forms the trivialisation of $E$ over all
of $\Rf$ (and hence $\sphere{4}$) except at the points where $x\in\Rf\iso\quat$
is a left quaternionic eigenvalue for $T$. However such points are merely
gauge singularities.\par
Since $E$ is an $\Su{2}$-bundle, we can simplify things a little here by
converting the data into data on an $\SP{1}$-bundle. We
can identify the fibre $\complexes^2$ with $\quat$, so that as complex vector spaces
\begin{eqnarray*}
\Hom{\complexes^2}{\complexes^k\tensor_\complexes\quat}&\iso&
\Hom{\quat}{\complexes^k\tensor_\complexes\quat}\\
&\iso& \complexes^k\tensor_\complexes\quat\tensor_\complexes\quat\\
&\iso& \complexes^k\tensor_\complexes\left(\quat\tensor_\reals\complexes\right)\\
&\iso& \complexes^k\tensor_\reals\quat\\
&\iso& \quat^k\tensor_\reals\complexes.
\end{eqnarray*}
This means that $\Hom{\complexes^2}{\complexes^k\tensor_\complexes\quat}$ may be regarded
as the complexification of the real space
\[
\left(\Rf\right)^k\iso\quat^k.
\]
So,using this conversion, we can recover ASD connections by choosing
\[
(T,P)\in\left(\Symm{\reals^k}{2}\tensor\quat\right)\oplus\quat^k=:\mr{k}
\]
subject to the ADHM condition
\begin{equation}\label{ADHMcond}
\Im\left(T^*T+PP^*\right)=0
\end{equation}
and the nondegeneracy condition, 
\begin{equation}\label{Nondeg}
\R_x \hbox{ is surjective for all } x\in\Rf.
\end{equation}
We regard $T$ as a symmetric matrix and $P$ as a column vector each with
quaternionic coefficients.
Now from a result by Wood \cite{WO}, every quaternionic matrix has a
left quaternionic eigenvalue. This means that $(T,0)$ always gives a 
reducible solution, i.e $(T,0)$ does not satisfy (\ref{Nondeg}). \par
We set $\Abad{k}\subset\mr{k}$ to be the
elements that satisfy (\ref{ADHMcond}) and define 
$\A{k}$ to be the set of all elements of $\Abad{k}$ that satisfy the 
nondegeneracy condition (\ref{Nondeg}).
\section{Group actions and Moduli}
It is a well known result that factoring out certain group actions from $\mc{k}$
or $\mr{k}$ we recover precisely the moduli space of instantons of charge
$k$. Essentially this says that an orbit of ADHM data with respect to a certain
group yields a gauge equivalence class of instantons.\par
For $\mr{k}$ this group is $\OSP{k}$ which acts on $\mr{k}$ via
\[
(\alpha,\beta):\ (T,P)\mapping(\alpha T\inv{\alpha},\alpha P\inv{\beta}).
\] 
It isn't hard to show that in fact
\[
\Im(T^*T+PP^*)
\]
is invariant with respect to the $\SP{1}$ action and equviariant with respect
to the $\O{k}$ action. In fact
\[
\vec{\mu}(T,P)=\Im(T^*T+PP^*)
\]
is the HyperK\"ahler moment map on $\mr{k}$ with its obvious hyperK\"ahler
structure with respect to the action of
$\O{k}$.\par
We fix the canonical basis of $\quat$ by
\begin{eqnarray*}
q_0&=&\id\\
q_1&=&\I\\
q_2&=&\J\\
q_3&=&\K.
\end{eqnarray*}
% \par
%Now the product Lie group $\OSP{k}$ acts on $\A{k}$ by
%\[
%(\alpha,\beta):\ (T,P)\mapping(\alpha T\inv{\alpha},\alpha P\inv{\beta}).
%\]
\begin{prop1}\label{freesp}
The action of $\SP{1}$ is free on $\A{k}$.
\end{prop1}
\begin{prf}
Suppose there is $\beta\in\SP{1}$ with $P\beta=P$. Since we require 
$(T,P)\in\A{k}$,we are assured that $P\neq0$. Thus 
$P^\mu\beta=P^\mu$ for some nonzero component $P^\mu\in\quat$ of $P$. Hence $\beta=1$.
\end{prf}
To show that $\O{k}$ acts freely on $\A{k}$, we need the following result.
%\begin{prop1}\label{SUk}
%Let $H$ be a compact matrix Lie Group and $\{\xi_i\}_{i=1}^l$ be 
%a set of 
%$l$ elements of $\goth{h}$. Suppose there 
%is $g\in H$ such that $Ad(g)\xi_i=\xi_i$ for all $i$, then there is a 
%splitting of $\goth{h}$ under which
%\[
%\xi_i=\left(\begin{array}{cc} \xi_i' & 0 \\ 0&0\end{array}\right),
%\]
%and 
%\[
%g=\left(\begin{array}{cc} g' & 0 \\ 0& g''\end{array}\right).
%\]
%\end{prop1}
%\begin{prf}
%Let $\goth{g}$ be the Lie subalgebra of $\goth{h}$ generated by 
%$\{\xi_i\}_{i=1}^l$ and $G$ the Lie subgroup of $H$ that it generates.
%We may use the Killing form to create the orthogonal complement 
%$\goth{g}^\perp$ of $\goth{g}$ in $\goth{h}$. We can show that 
%$\goth{h}=\goth{g}\oplus\goth{g}^\perp$ is a direct sum of Lie algebras. 
%Hence $\goth{g}^\perp$ 
%generates a Lie subgroup $G^\perp$ which commutes with $G$, that is acts 
%trivially on $\goth{g}$. Hence $\xi_i\in\goth{g}$ has the required form,
%and any $g$ which fixes that space must also have the prescribed form.
%\end{prf}
\begin{lem1}\label{reducibles}
If $(T,P)\in\mr{k}$
is fixed by $u\in\O{k}$ then there is a decomposition
\[
T=\left(\begin{array}{cc} T' &0\\ 0 & T''\end{array}\right)\quad 
P=\left(\begin{array}{c}  P' \\ 0 \end{array}\right)
\]
with
$(T',P')\in\mr{l}$ and $T''\in\Symm{\reals^{k-l}}{2}\tensor\quat$.
%If the fixed set of $(T,P)$ contains two or more distinct elements then 
%${T''}^*T''$ is a real matrix. 
Hence such $(T,P)\notin\A{k}$.
\end{lem1}
\begin{prf}
Suppose $u\in\O{k}$ fixes $(T,P)$.
%We consider first the action on the $T$ part. We may write
%\[
%T={T_\goth{o}}+\frac{1}{k}\hbox{tr}T\id
%\]
%where $T_\goth{o}\in(\Symm{\reals^k}{2})_0\tensor_\reals\quat$, the traceless
%symmetric endomorphisms. Since $\O{k}$ preserves the trace part, we need only
%look at the the traceless part of $T$
%We can make the identification
%\[
%\Skew{\reals^k}{2}\mapping(\Symm{\reals^k}{2})_0\cr
%v\wedge w\mapsto v\odot w
%\]
%which preserves the action of $\O{k}$. Note that 
%$\Skew{\reals^k}{2}\iso\o{k}$.

%We need to look at  the ``complex case'', that is 
%the space  $\left(i\u{k}\oplus\complexes^k\right)\tensor_\reals\quat$ acted on 
%by the group $\U{k}\times\SP{1}$. 
%By our previous argument, this will reduce to the 
%result stated.\par
%Let $(T,P)$ be in this ``complexified'' space 
%$\left(i\u{k}\oplus\complexes^k\right)\tensor_\reals\quat$ and
%suppose that there is $\alpha\in\U{k}$ such that
%\[
%\alpha T\inv{\alpha}=T,\quad \alpha P=P.
%\]
Then $u P=P$, and $P={P_i}{q_i}$ for vectors $P_i\in\reals^k$, and
$u P_i= P_i$ for each $i=0...3$.
Hence $u$ has at least one eigenvector with eigenvalue 1. 
Decompose $\reals^k$ into $V\oplus W$ where $V$ is the maximal $1$-eigenspace 
of $u$ and $W$ its orthogonal complement. Notice that $V\supset\mbox{Span}\{P_i|i=0...3\}$. 
With respect to this decomposition, we have
\[
T=\left(
\begin{array}{cc}
T'&T_0^\top\\
T_0&T''
\end{array}
\right),\quad
P=\left(
\begin{array}{c}
P'\\
0
\end{array}
\right),
\]
and
\[
u=\left(
\begin{array}{cc}
\id & 0\\
0 & u'
\end{array}
\right),
\]
for some $u'\in\O{k-l}$.
Assume that $u\neq\id$ and hence $W\neq 0$. The condition that $uT\inv{u}=T$ 
shows us, in 
particular that $u'T_0=T_0$, so the columns of $T_0$, $T_0^\mu$ say, satisfy
\[
u'T_0^\mu=T_0^\mu
\]
But $u'$ does not have $+1$ as an eigenvalue by the decomposition. Hence
$T_0^\mu=0$ and $T_0=0$.\par
Thus we have the following decompositions for $T$
and $P$ with respect to this splitting.
\[
T=\left(\begin{array}{cc} T' &0\\ 0 & T''\end{array}\right)\quad 
P=\left(\begin{array}{c}  P' \\ 0 \end{array}\right).
\]
Now by R. Wood \cite{WO}, we know that $T''$ (hence $T$) has a left eigenvalue
$\lambda\in\quat$. So 
\[
\F_{(\lambda,T,P)}=\R_{(\lambda,T,P)}\R_{(\lambda,T,P)}^*
\]
will not be invertible and thus $(T,P)\notin\A{k}$.

%Thus $T$ and $P$ have the required form.
\end{prf}
\begin{cor1}
The following are equivalent
\begin{enumerate} 
\item $(T,P)$ satisfies the nondegeneracy condition;
\item $(T,P)$ has trivial stabiliser under $\O{k}$.
\end{enumerate}
\end{cor1}
\begin{cor1}\label{free}
The action of $\O{k}$ is free on $\A{k}$.
\end{cor1}
We may have a small problem here. Although $\O{k}$ and $\SP{1}$ individually act freely on $\A{k}$, the full group $\O{k}\times\SP{1}$ doesn't act freely here. Each point is fixed by $\pm(\id,1)$ meaning that the group of symmetries we require is 
\[
\frac{\O{k}\times\SP{1}}{\integers_2}.
\]
%Of course, if $k$ is odd we have
%\[
%\O{k}=\So{k}\times\integers_2
%\]
%hence,
%\[
%\frac{\O{k}\times\SP{1}}{\integers_2}\iso\So{k}\times\SP{1}.
%\]
Not only that but for $k\ge2$, we have some $(T,P)\in\mr{k}
%\left(\Symm{\reals^k}{2}\oplus\reals^k\right)\tensor\quat
$
such that 
\[
gT\inv{g}=T,\quad gP\inv{q}=P
\]
for certain $(g,q)\in\OSP{k}$.
 We must check that they
are not in $\A{k}$.
\subsection{The Action of $\OSP{k}$ on $\mr{k}$}
Choose $\xi\in\o{k}$, and choose a basis of $\reals^k$ (hence of $\quat^k$) in 
which 
\[
\xi=\left(
\begin{array}{ccc}
0    & \id_l & 0 \\
-\id_l & 0   & 0 \\
0 & 0& 0_{k-2l} 
\end{array}
\right)
\]
where $\id_l$ is the $l\times l$ identity matrix and $0_{k-2l}$ is the $(k-2l)\times
(k-2l)$ zero matrix. It may happen that $k=2l$, in which case $\xi$ is invertible
and forms a complex structure on $\reals^k=\reals^{2l}$.
Thus for
\[
T=\left(
\begin{array}{ccc}
 T' & T_1 & T_2^\top \\ 
T_1^\top & T'' & T_3^\top \\
T_2& T_3 & T'''
\end{array}
\right),\quad
P=\left( 
\begin{array}{c}
P' \\
P''\\
P'''\\
\end{array}
\right)
\]
we find
\[
\left[\xi,T\right]=\left(
\begin{array}{ccc}
 T_1+T_1^\top & T''-T' & T_3^\top \\ 
T'-T''&-T_1-T_1^\top & -T_2^\top \\
T_3& -T_2 & 0
\end{array}
\right),\quad
\xi P=\left( 
\begin{array}{c}
-P'' \\
P'\\
0\\
\end{array}
\right)
\]
So if $[\xi,T]=0$ and $\xi P=P\alpha$ then 
\[
T=\left(
\begin{array}{ccc}
 T' & T_1 & 0 \\ 
-T_1& T'  & 0 \\
   0& 0   & T'''
\end{array}
\right),\quad
P=\left( 
\begin{array}{c}
P' \\
P'\alpha\\
0\\
\end{array}
\right)
\]
with the condition that $|\alpha|=1$.
From this we can conclude, that if $(T,P)$ satisfies the ADHM condition then it
is a reducible solution, unless $\xi$ is invertible and hence $k$ is even. \par
Let us now suppose that $\xi$ is invertible, then as we mentioned, $\xi$
is a complex structure on $\reals^{2l}$.
Let us take the complex point of view and look at the process for 
$\complexes^l$. Under this identification $\xi$ becomes multiplication by $i$ 
and
\[
T=T'+iT'',P=P'+iP''.
\]
Also our condition that $\xi P=P\alpha$ becomes
\[
iP=P\alpha
\]
or
\[
-P'+iP'=P'\alpha+iP'\alpha.
\]
Comparing the real and imaginary parts,
\begin{eqnarray*}
P'\alpha&=&-P',\\
P'\alpha&=&P',
\end{eqnarray*}
since $\alpha$ is considered to have only real coefficients.
Hence $P'=0$ and $P=0$, and the solution is again reducible.
%that in our chosen basis of $\reals^k$, we have 
%\[
%\xi=\left(
%\begin{array}{cc}
%0    & \id  \\
%-\id & 0   \\ 
%\end{array}
%\right),\quad
%T=\left(
%\begin{array}{cc}
% T' & T'' \\ 
%-T''& T'  \\
%\end{array}
%\right),\quad
%P=\left( 
%\begin{array}{c}
%P' \\
%P'\alpha\\
%\end{array}
%\right).
%\]
%Thus
%\[
%\R_{0}\R_{0}^*=\left(
%\begin{array}{cc}
%(T')^*(T')-(T'')^*(T'')+(P')(P')^* & (T')^*(T'')+(T'')^*(T')-(P')\alpha(P')^*\\
%-(T'')^*(T')-(T')^*(T'')+(P')\alpha(P')^*& -(T'')^*(T'')+(T')^*(T')+(P')(P')^*
%\end{array}
%\right).
%\]
%If $(T,P)\in\A{k}$ then 
%\[
%\Im\left((T')^*(T'')+(T'')^*(T')-(P')\alpha(P')^*\right)=0
%\]
%Let us examine any diagonal element of the left hand side:
%\[
%\Im\left(\sum_{j=0}^k \left((T')_{ij}^*(T''){ji}+(T''){ij}^*(T')_{ji}\right)-(P')_i\alpha(P')_i^*\right).
%\]
%Since $T'$ is symmetric  and $T''$ is skewsymmetric, we must have
%\begin{eqnarray*}
%& &\Im\left(\sum_{j=0}^k \left((T')_{ij}^*(T''){ji}+(T''){ij}^*(T')_{ji}\right)-(P')_i\alpha(P')_i^*\right)\\
%&=&\Im\left(\sum_{j=0}^k \left((T')_{ji}^*(T''){ji}+(T''){ij}^*(T')_{ji}\right)-(P')_i\alpha(P')_i^*\right)\\

%\[
%\frac{\As{k}}{\O{k}}\bigg/\SP{1}=\frac{\As{k}}{\OSP{k}}\bigg/\integers_2.
%\]
%Our method here is to perform the factorisation of $\A{k}$ by $\O{k}$ giving
%rise to the framed moduli space and then to factorise by $\SP{1}$. 
%\end{rmk}
So far, we have proved that if the vector field induced by 
$(\xi,\alpha)\in\o{k}\oplus\sP{1}$ vanishes at a point $(T,P)$ then
$(T,P)$ is a reducible solution. This in turn shows that the stabiliser of
$(T,P)\in\A{k}$ must be a discrete, hence finite subgroup of $\OSP{k}$ which
is enough for our purposes.
\subsection{Introducing the Moduli Spaces}
Let 
\[
\tilde{M}_{k}=\frac{\A{k}}{\O{k}}
\]
and  
\begin{eqnarray*}
M_k&=&\tilde{M}_{k}/\left(\frac{\SP{1}}{\integers_2}\right)\\
   &=&\tilde{M}_{k}/\So{3}.
\end{eqnarray*}
Define the Atiyah map $\function{\At}{\A{k}}{\mathcal{A}_k}$ where 
$\mathcal{A}_k$ is the space of connections of charge $k$ by
\[
(T,P)\mapsto v^*\hbox{d}v
\]
as above.
It is well known that any two elements of the same orbit produce gauge 
equivalent connections, and that $M_k$ is diffeomorphic to $\mathcal{M}_k$, the
moduli space of ASD connections of charge $k$\cite{N}.\par
Also, by considering $\BBt{k}=\frac{\A{k}}{\O{K}}$, we obtain the framed moduli
space of connections. The manifold of equivalence classes of ADHM data under 
the action of $\O{k}$ is precisely the moduli space of framed instantons,
$\widetilde{\mathcal{M}}_k$.\par
%Later on, we will often write $\mathcal{M}_k$ and $\widetilde{\mathcal{M}}_k$
%for $M_k$ and $\tilde{M}_k$ when we have show that that the two spaces give
%rise to the same bundle data.
\subsection{Curvature of bundles under group actions}\label{Equivbun}
Here, we recall the theory of vector bundles and group actions detailed in 
section 5.2.3 of \cite{DK}. Let $\function{\widehat{\pi}}{\hE}{\widehat{Y}}$ be
a vector bundle, and $G$ a Lie group whose action on $\tE$ covers a free action
on $\widehat{Y}$. Also let $\hE$ be endowed with a $G$-invariant connection
$\hcon$. Our aim is to construct a connection $\grad$ on the factor bundle 
$E=\frac{\hE}{G}$ over $Y=\frac{\widehat{Y}}{G}$.\par
To do this we need a connection on the principal $G$-bundle 
$p:\widehat{Y}\mapping Y$. This will enable us to lift tangent vectors on $Y$
to $\widehat{Y}$ and compute directional derivatives. We suppose that it is
given in the form of a horizontal distribution $H\subset\Tang{\widehat{Y}}$,
with connection 1-form $\theta$. \par
Any section $s\in\sect{Y}{E}$ comes from an invariant section 
$\widehat{s}\in\sect{\widehat{Y}}{\hE}$. Therefore we can set
\[
\widehat{\grad_X s}=\hcon_{\widehat{X}}\widehat{s},
\]
where $\widehat{X}$ is the horizontal lift of $X$. This descends to the quotient
so that $\widehat{\grad_X s}$ is the lift of the object which can be called
$\grad_X s$.\par
From the definition of $Y$, we know that $p^*E\iso\hE$, so we can consider the 
effect of the pull 
back $p^*\grad$ on $\hE$. Now, the directional derivatives of $p^*\grad$ will
vanish on vertical vectors in $\Tang{\widehat{Y}}$, hence we know that
\[
\hcon=p^*\grad+V
\]
where $V$ vanishes on horizontal vectors, and is equivariant under the $G$ 
action by definition of $\hcon$.\par
It is therefore obvious that $V=\Psi\theta$ where $\Psi$ is a linear
$\mathfrak{g}\mapping\End{\hE}$\par
It can then be shown (and \cite{DK} do this to some extent) that
\begin{equation}\label{QF}
F(\hcon)(\widehat{X_1},\widehat{X_2})=F(\grad)(X_1,X_2)+\Psi F(\theta)(\widehat{X_1},\widehat{X_2}),
\end{equation}
allowing us to compare the curvature of $\hcon$ with $\grad$.
\section{The case of charge 1 instantons ($k=1$)}\label{k=1}
For $k=1$, we are in a truly interesting position since for any
\[
(T,P)\in\A{1}\subset\mr{1}=\left(\Symm{\reals}{2}\tensor\quat\right)\oplus\quat\iso\quat\oplus\quat,
\]
the ADHM condition
\[
\Im(T^*T+PP^*)=0
\]
is automatically satisfied, hence $\Abad{1}=\mr{1}
%\Symm{\reals}{2}\tensor\quat\oplus\quat
$.\par

Now we can construct a canonical bundle $\hE\mapping\As{1}\times\Rf$ as
follows. \\
For each $(T,P,x)\in\As{1}\times\Rf$, define the fibre of $\hE$ to be
\[
\hE_{(T,P,x)}=\ker\R_x=\ker(T-x,P):\quat\oplus\quat\mapping\quat.
\]
We can define a connection $\hcon$ on $\hE$ given at the point $(T,P,x)$ by
$\At(T,P)_x$. From the suggestive terminology, it is clear that we will choose
$\widehat{Y}=\A{1}\times\Rf$. Since we have removed the singular points, the action of $$\frac{\OSP{1}}{\integers_2}=\frac{\integers_2\times\SP{1}}{\integers_2}=\SP{1}$$ is free. It is also clear that by construction
$\hcon$ is $\SP{1}$ invariant.\par
%Also $$\tilde{M}^*_{1}=\frac{\A{1}}{\O{1}}=\frac{\A{1}}{\integers_2},$$
%where $\integers_2$ acts via
%\[
%\rho:(T,P)\mapsto(\rho T\inv{\rho},\rho P)=(T,\rho{P})
%\]
%So, if $(T,P)\in\A{1}$, then $(T,[P])\in\tilde{M}^*_{1}$ where
%\[
%[P]=[-P].
%\]
\subsection{The curvature of the $\SP{1}$-bundle}
We now need to consider the principal $\SP{1}$-bundle
\[
\A{1}\mapping\BBs{1}.
\]
For $(T,P)\in\A{1}$, the action of the group $\SP{1}$ is
\[
u:(T,P)\mapsto(T,P\inv{u}).
\]
The vertical subspace $V_{(T,P)}$ will therefore be
\[
\{(0,-P\widehat{u})\vert\widehat{u}\in\sP{1}\}.
\]
We define a horizontal subspace $H_{(T,P)}$ to be the orthogonal complement of
the vertical subspace. Thus $H_{(T,P)}$ will be the space
\[
\{(t,p)\vert P^*p\in\reals\}=\ker(p\mapsto\Im(P^*p))\}.
\]
This immediately gives us the connection 1-form
\begin{equation}\label{Aconn}
\theta_{(T,P)}(t,p)=-\frac{1}{\mod{P}^2}\Im(P^*p).
\end{equation}
\begin{prop1} $\theta$ is indeed a connection 1-form.
\end{prop1}
\begin{prf}
First, $\theta$ is well defined as $P\neq 0$.\par
On vertical vectors
\begin{eqnarray*}
\theta_{(T,P)}(0,-P\widehat{u})&=&-\frac{1}{\mod{P}^2}\Im(-P^*P\widehat{u})\\
                               &=&\Im(\widehat{u})\\
                               &=&\widehat{u}.
\end{eqnarray*}
Also, for equivariance,
\begin{eqnarray*}
\theta_{(T,P\inv{u})}(t,p\inv{u})&=&-\frac{1}{\mod{Pu}^2}\Im(uP^*p\inv{u})\\
                           &=&-\frac{1}{\mod{P}^2}u\Im(P^*p)\inv{u}\\
                           &=& \hbox{ad}(u)\theta_{(T,P)}(t,p).
\end{eqnarray*}
\end{prf}
Using local coordinates $(T,P)$, this connection 1-form can be written
\begin{eqnarray*}
\theta=-\frac{1}{\mod{P}^2}\Im(P^*\hbox{d}P),
\end{eqnarray*}
so the curvature form restricted to the horizontal space is given by
\begin{eqnarray*}
F(\theta)\restrict{\ker\theta}&=&\d\theta\restrict{\ker\theta}\\
                              &=&-\frac{1}{\mod{P}^2}\d P^*\wedge\d P\restrict{\ker\theta}
-\frac{1}{\mod{P}^2}\d\mod{P}^2\wedge\theta\restrict{\ker\theta}\\
&=&-\frac{1}{\mod{P}^2}\d P^*\wedge\d P\restrict{\ker\theta}.
%&=&\frac{\hbox{d}\mod{P}^2}{\mod{P}^4}\wedge P^*\hbox{d}P-\frac{1}
%{\mod{P}^2}\hbox{d}P^*\wedge\hbox{d}P-\frac{1}{\mod{P}^4}P^*\hbox{d}P\wedge
%P^*\hbox{d}P.
\end{eqnarray*}
%Now we have to consider the behaviour of $F(\theta)$ on the horizontal
%subspace. The point to realise here is that the horizontal subspace is 5
%dimensional, with 4 dimensions in the $\Symm{\reals^1}{2}\tensor\quat$
%direction and the remaining 1 dimension in the $\quat^1$ direction. So horizontal
%tangent vectors are one dimensional in the $P$-direction. The curvature
%form however is independent of the $T$ directions and a two form in the $P$
%directions which is one dimensional hence the two form vanishes when restricted
%to the Horizontal subspace.\par
Now $(t,p)\in\Abad{1}$ is horizontal if and only if $p=\lambda P$
by the definition of $\theta$ in (\ref{Aconn}). Hence for a
horizontal vector $(X,t,\lambda P)\in\Rf\times H_{(T,P)}$
\[
\d P((X,t,\lambda P))=\lambda P,
\]
and from this we can see that
\begin{equation}\label{dP}
\d P\restrict{\ker\theta}=\frac{\delta\mod{P}^2}{2\mod{P}^2}P=\frac{\delta\mod{P}}{\mod{P}}\restrict{\ker\theta}P.
\end{equation}
As a result
\[
F(\theta)\restrict{\ker\theta}=-\delta\mod{P}\wedge\delta{\mod{P}}\restrict{\ker\theta}=0.
\]
This means that the formula (\ref{QF}) becomes
\[
F(\con)(X_1,X_2)=F(\hcon)(\widehat{X_1},\widehat{X_2}).
\]
\subsection{The Curvature of the Universal Bundle}\label{1curv}
Recall that the Atiyah map was defined by
\begin{eqnarray*}
&\function{\At}{\As{1}}{\mathcal{A}_k}\\
&(T,P)\mapsto v(T,P)^*\hbox{d}v(T,P).
\end{eqnarray*}
Explicitly for $k=1$, $\At(T,P)$ is given by $v(T,P)^*\hbox{d}v(T,P)$ where
\[
v(T,P)=
\frac{\mod{x-T}}{\sqrt{\mod{P}^2+\mod{x-T}^2}}
\left[
\begin{array}{c}
\frac{(x-T)}{\mod{x-T}^2}P\\
1
\end{array}
\right].
\]
Thus $\hcon$ is given at the point $(T,P,x)$ by $v(T,P)_x^*\widehat{\hbox{d}}
v(T,P)_x$, where
$\widehat{\hbox{d}}$ denotes the de Rham differential not only with respect to
$x$ but also with respect to the quaternions $T$ and $P$. We will also denote by $\delta$ the de Rham differential in the ADHM space. Also, set $\tilde{x}=x-T$
\par
Expanding this calculation and setting
\[
\Delta={\sqrt{\mod{P}^2+\mod{\tilde{x}}^2}},\]
we find that the connection 1-form of $\hcon$ with respect to the global
trivialisation of $\hE$ is
\begin{equation}\label{udu}
v(T,P)_x^*\widehat{\hbox{d}} v(T,P)_x=
\frac{1}{\Delta^2}\Im\left(P^*\d
P-\frac{1}{\mod{\tilde{x}}^2}P^*\hd\tilde{x}^*\tilde{x}P\right).
\end{equation}
We compute the curvature in the obvious way, i.e.
\[
F(A)=\d A+A\wedge A.
\]
From this we get
\begin{equation}\label{Fhat}
\begin{array}{c}
F(\hcon)=\frac{1}{\Delta^4}\bigg[ \mod{\tilde{x}}^2\d P^*\wedge\d
P +\frac{1}{\mod{\tilde{x}}^2}P^*
\tilde{x}^*\hd\tilde{x}\wedge\hd\tilde{x}^*\tilde{x}P\\
 -
P^*\tilde{x}^*\hd\tilde{x}\wedge\d P-\d
P^*\wedge\hd\tilde{x}^*\tilde{x} P \bigg].
\end{array}
\end{equation}
Now, we are interested in ``factoring'' out $\hcon$ by $\SP{1}$. That is, we
want only to consider the horizontal lifts of tangent vectors and vector fields
on $\BBs{1}$ and how $F(\hcon)$ behaves when restricted to these.\par
%Now $(t,p)\in\Abad{1}$ is horizontal if and only if $p=\lambda P$ by the
%definition of $\theta$ in (\ref{Aconn}).Hence for a horizontal vector
%$(X,t,\lambda P)\in\Rf\times H_{(T,P)}$
%\[
%\d P((X,t,\lambda P))=\lambda P,
%\]
%and from this we can see that
%\[
%\d P\restrict{\ker\theta}=\frac{\delta\mod{P}^2}{2\mod{P}^2}P=\frac{\delta\mod{P}}{\mod{P}}P
%\]
So by (\ref{dP}), we have
\begin{eqnarray*}
& &F(\hcon)\restrict{\ker\theta}\\
&=&\frac{1}{\Delta^4}\bigg[\mod{\tilde{x}}^2{\delta\mod{P}\wedge\delta\mod{P}}
+\frac{1}{\mod{\tilde{x}^2}}P^*\tilde{x}^*\hd\tilde{x}\wedge\hd\tilde{x}^*\tilde{x}P
-2\Im(P\tilde{x}^*\hd\tilde{x}\wedge\frac{\delta\mod{P}}{\mod{P}}P) \bigg]\\
&=&\frac{1}{\Delta^4}\bigg[\frac{1}{\mod{\tilde{x}^2}}P^*\tilde{x}^*\hd\tilde{x}\wedge\hd\tilde{x}^*\tilde{x}P
+2\frac{\delta\mod{P}}{\mod{P}}\wedge\Im(P\tilde{x}^*\hd\tilde{x}P) \bigg].\\
\end{eqnarray*}
If we set $\mod{P}=\rho$, then we find that when restricted to the horizontal
space
\begin{equation}\label{Aconn2}
F(\hcon)=\frac{1}{\Delta^4}\bigg[\frac{1}{\mod{\tilde{x}^2}}P^*\tilde{x}^*\hd\tilde{x}\wedge\hd\tilde{x}^*\tilde{x}P
+2\frac{\delta\rho}{\rho}\wedge\Im(P\tilde{x}^*\hd\tilde{x}P)
\bigg].
\end{equation}
This is very similar to Anselmi's formula in \cite{An}. In fact if we assume,
(as he effectively does) that $P$ is real and positive under part of a gauge
fixing condition, then we have his answer scaled by a factor of $\frac{1}{2}$.
\Endoff{Section2.2.tex}
\subsection{Calculating the $\mu$ map}\label{1mu}
Our next goal is to compute a form representing 
the second Chern class of $E$. 
This will be the same as the second Chern character $\hbox{ch}_2(E)$ since
the group is $\Su{2}\iso\SP{1}$ and hence $\hbox{c}_1(E)=0$.
We may obtain a representative of the cohomology class of $c_2(E)$ given by
\[
\frac{1}{4\pi^2}\hbox{tr}\left(F(\grad)\wedge F(\grad)\right).
\]
By expanding this we find that the desired form is
\begin{equation}\label{c2e}
c=\frac{6\rho^4}{\Delta^8\pi^2}\hd\tilde{x}_1\wedge\hd\tilde{x}_2\wedge\hd\tilde{x}_3\wedge\hd\tilde{x}_4-\frac{6\rho^3}{\Delta^8\pi^2}\delta\rho\wedge\sum_{i=1}^{4}(-1)^i\tilde{x}_i\hd\tilde{x}_1\ldots\wedge\hat{i}\wedge\ldots\wedge\hd\tilde{x}_4
\end{equation}
where $\tilde{x}_i=x_i-T_i$.\par
Now, it can be checked that 
\[
c=\hd\left(\frac{1}{2\pi^2}\frac{\left(|\tilde{x}|^2+3\rho^2\right)}{\left(|\tilde{x}|^2+\rho^2\right)^3}\sum_{i=1}^{4}(-1)^i\tilde{x}_i\hd\tilde{x}_1\ldots\wedge\hat{i}\wedge\ldots\wedge\hd\tilde{x}_4\right).
\]
We choose a compact $d$-dimensional submanifold $\Sigma$ of $\reals^4$ and let $\alpha\in\Omega_{cpt}^{4-d}(\reals^4)$ be a closed form such that for
any $\beta\in\Form{\Sigma}{d}$. 
\[
\int_\Sigma\beta=\int_{\Rf}\beta\wedge\alpha.
\] 
We also choose a submanifold $\Xi$ of $M_1$ whose intersection with $\partial
M_1$ is compact.
\begin{eqnarray*}
\int_{\Xi\times\Sigma} c
&=&\int_{\Xi\times\Rf}c\wedge\alpha\\
&=&\int_{(\partial M\cap\Xi)\times\Rf}\left(\frac{1}{2\pi^2}\frac{\left(|\tilde{x}|^2+3\rho^2\right)}{\left(|\tilde{x}|^2+\rho^2\right)^3}\sum_{i=1}^{4}(-1)^i\tilde{x}_i\hd\tilde{x}_1\ldots\wedge\hat{i}\wedge\ldots\wedge\hd\tilde{x}_4\right)\wedge\alpha\\
&=&\frac{1}{2\pi^2}\int_{(\Xi\cap\{\rho=0\})\times\Rf}\frac{1}{|\tilde{x}|^4}\sum_{i=1}^{4}(-1)^i\tilde{x}_i\hd\tilde{x}_1\ldots\wedge\hat{i}\wedge\ldots\wedge\hd\tilde{x}_4\wedge\alpha\\
&=&\frac{1}{2\pi^2}\int_{(\Xi\cap\{\rho=0\})\times\Sigma}\frac{1}{|\tilde{x}|^4}\sum_{i=1}^{4}(-1)^i\tilde{x}_i\hd\tilde{x}_1\ldots\wedge\hat{i}\wedge\ldots\wedge\hd\tilde{x}_4.
\end{eqnarray*}
The right hand side above can be seen to be the Gau{\ss} formula for the 
linking number of $\Sigma$ with $\Xi\cap\{P=0\}$ regarded as a 3-dimensional 
submanifold of $\Rf$.
\subsection{The Donaldson Polynomials in the case $k=1$}
We now have to assess the consequences of this for the Donaldson Polynomials.
Let $\Sigma_1,\ldots,\Sigma_d$ be $d$ submanifolds of $\Rf$. We may form the
Donaldson $\mu$ class 
\[
\mu_{Don}(\Sigma_i)=\int_{\Sigma_i}c
\]
and from this the Donaldson polynomial
\begin{eqnarray*}
\mathrm{Don}_1(\Sigma_1,\ldots,\Sigma_d)
&=&\int_{\mathcal{M}_1}\mu_{Don}(\Sigma_1)\ldots\mu_{Don}(\Sigma_d)\\
&=&\int_{\mathcal{M}_1\times\Sigma_1\times\ldots\times\Sigma_d}
c_1\wedge\ldots\wedge c_d
\end{eqnarray*}
where $c_i$ is $c$ restricted to $\Sigma_i$.\par
Now let
\[
f(x,T,\rho)=\frac{1}{2\pi^2}\frac{\left(|\tilde{x}|^2+3\rho^2\right)}{\left(|\tilde{x}|^2+\rho^2\right)^3}
\]
so that 
\[
c=\hd\left(f(x,T,P)\alpha(x,T,P)\right)
\]
where
\[
\alpha(x,T,P)=\sum_{i=1}^{4}(-1)^i\tilde{x}_i\hd\tilde{x}_1\ldots\wedge\hat{i}\wedge\ldots\wedge\hd\tilde{x}_4.
\]
Then
\begin{eqnarray*}
& &\mathrm{Don}_1(\Sigma_1,\ldots,\Sigma_d)\\
&=&\int_{\mathcal{M}_1\times\Sigma_1\times\ldots\times\Sigma_d}
c_1\wedge\ldots\wedge c_d
\\
&=&\int_{\mathcal{M}_1\times\Sigma_1\times\ldots\times\Sigma_d}
\hd(f(x^1,T,P)\alpha(x^1,T,P))\wedge\ldots\wedge\hd(f(x^d,T,P)\alpha(x^d,T,P))\\
&=&\int_{\mathcal{M}_1\times\Sigma_1\times\ldots\times\Sigma_d}
\hd\left(f(x^1,T,P)\alpha(x^1,T,P)\wedge\ldots\wedge\hd(f(x^d,T,P)\alpha(x^d,T,P))\right)\\
&=&\int_{\{\rho=0\}\times\Sigma_1\times\ldots\times\Sigma_d}
f(x^1,T,P)\alpha(x^1,T,P)\wedge\ldots\wedge\hd(f(x^d,T,P)\alpha(x^d,T,P))\\
&=&\frac{1}{2\pi^2}\lim_{\rho\mapping
0}\int_{\{(T,\rho)\}\times\Sigma_1\times\ldots\times\Sigma_d}\frac{(3\rho^2+|x^1-T|^2)}
{(\rho^2+|x^1-T|^2)^3}\alpha(x^1,T,0)\wedge\\
& &\prod_{l=1}^d\left(\frac{6\rho^4}{(\rho^2+|x^l-T|^2)^4}
\hd(x_1^l-T_1)\wedge\ldots\hd(x_4^l-T_4)\right)
\end{eqnarray*}
using the formula (\ref{c2e}) for $c$.\par
Now, since we have a singularity when $x=T$ and $\rho=0$ in the formula
\[
\frac{6\rho^4}{(\rho^2+|x-T|^2)^4}
\]
we have to be a bit careful with limits.
Now, a calculation shows that
\[
\lim_{\rho\mapping
0}\frac{6\rho^4}{(\rho^2+|x-T|^2)^4}=\frac{1}{2}\mathrm{vol}\,\sphere{3}\delta(x-T)
\]
and provided $x^1,\ldots,x^d$ are distinct points in $\Rf$
\[
\lim_{\rho\mapping
0}\prod_{l=1}^{d}\frac{6\rho^4}{(\rho^2+|x^l-T|^2)^4}=\prod_{l=1}^{d}\frac{1}{2}\mathrm{vol}\,\sphere{3}\delta(x^l-T).
\]
For $d=2$, we have
\[
\mathrm{Don}_1(\Sigma_1,\Sigma_2)=\frac{\mathrm{vol}\,\sphere{3}}{2}
\int_{\Sigma_1}\int_{\Sigma_2}\frac{1}{|x^1-x^2|^2}\alpha(x^1,x^2,0)
\]
which is a constant multiple of the linking number of $\Sigma_1$ with
$\Sigma_2$.\par
For $d>2$, the situation is much more complicated. A discussion of this can be found in another article by Anselmi, \cite{An2}.\par
%It is obvious that this discovery by Anselmi is anomalous. In classical %Donaldson
%theory, the Donaldson polynomials are polynomials in $\mathrm{H}^2(M;\reals)$,
%and since $\mathrm{H}^2(\sphere{4};\reals)=\{0\}$, all Donaldson polynomials
%must be trivial, contradicting Anselmi's results.
Although Donaldson theory is relevant only for 4-manifolds with $b_+^2>0$,
the formation of $\mu$-classes is completely trivial on $\sphere{4}$, and one would expect all polynomials constructed with Chern-Weil representatives of these $\mu$-classes as described above to agree with this triviality.
However, Anselmi has produced ``representatives" of $\mu$-classes which form
non-trivial polynomials, contradicting our expectations.
\Endoff{Section2.3.tex}
\subsection{Reasons for the Anomalies}
Let us examine what is going on more closely. The bundle $\hE$ is defined on
$\Rf\times\mc{k}$ minus the set $S_{k}$ consisting of the points $(x,T,P)\in\Rf\times\mc{k}$ such that there is $u\in\U{k}$ for which
\[
uT\inv{u}=\left(\begin{array}{cc}
T' &0 \\
0  &x
\end{array}\right),\quad
uP=\left(\begin{array}{c}
P'\\
0
\end{array}\right),
\]
where $(T',P')\in\mc{k-1}$. We are then integrating a form on the
quotient of the hyperK\"ahler reduction of
$\left(\Rf\times\mc{k}\right)\backslash S_k$, i.e.
\[
\mathbb{M}_k=\frac{\left((\Rf\times\mc{k})\backslash S_k\right)/\!\!/\!\!/\!\!/\U{k}}{\So{3}}.
\]
%(We will deal with hyperK\"ahler geometry in the next chapter.)
(Here $\So{3}$ acts trivially on the manifold $\Rf$.)
\par
In the infinite dimensional case, for a submanifold $\Sigma$ of $\Rf$
we form
\[
\mu_{Don}(\Sigma)=c_2(\mathbb{E})/[\Sigma]=\int_\Sigma c_2(\mathbb{E})
=\int_{\Rf} c_2(\mathbb{E})\wedge\pi^*\mathrm{PD}(\Sigma_i)\in\mathrm{H}^{4-\dim{\Sigma}}(\mathcal{M}^*)
\]
where $c_2(\mathbb{E})\in \mathrm{H}^4(\Rf\times\mathcal{M}^*)$
and $\pi: \Rf\times\mathcal{M}\mapping\Rf$ is the projection. The
slant product is well-defined as an integration along the fibre because we have a perfectly decent
trivial fibration $\Rf\times\mathcal{M}^*\mapping\Rf$.
\par However in the finite case the fibration
$$\mathbb{M}_k\mapping\Rf$$
is destroyed because we remove points
(namely $S_k$) from the direct product before taking the various quotients;
the preimage varies topologically from point to point. If we na\"ively form the finite
dimensional version
\[
\mu_{Don}(\Sigma)=\frac{-1}{8\pi^2}\int_\Sigma\mathrm{tr}(F(\grad)\wedge F(\grad))
\]
then it is not altogether clear where this $\mu_{Don}(\Sigma)$ lies.
%The space
% $\mathcal{M}$ isn't smooth and hence
The Chern-Weil representative of $\mu_{Don}(\Sigma)$  doesn't really
represent a cohomology class on $\mathcal{M}$. We are forced therefore to
reinterpret the situation for the ADHM case.\par
We have the 4-form $c$ representing the second Chern class of the universal
bundle $\hE$ over the manifold $\mathbb{M}_k$.
Let $\Sigma_1,\ldots,\Sigma_k$ be compact submanifolds of $\Rf$ without boundary whose dimensions sum
to $4l-8k+3$, that is of the correct dimensions to form a Donaldson polynomial.
We can then consider the form
\[
\mu_i=c\wedge\iota^*\pi^*\mathrm{PD}(\Sigma_i)
\]
where $\pi:\Rf\times\mathcal{M}_{k}\mapping\Rf$ is projection, and
$\iota:\mathbb{M}_k\hookrightarrow\Rf\times\mathcal{M}_k$
inclusion. Thus we can form the integral
\[
\int_{D_{(k,l)}}\Delta^*\left(\Pi_1^*\mu_1\wedge\ldots\wedge\Pi_l^*\mu_l\right)
\]
where $\Pi_i:\mathbb{M}_k^l\mapping \mathbb{M}_k$ is projection onto the $i$th factor,
\[
D_{(k,l)}=\{(x_1,\ldots,x_l,[T,P])\in(\Rf)^l\times\mathcal{M}_{k}|(x_i,[T,P])\in \mathbb{M}_k\hbox{ for all }k\}
\]
and
\begin{eqnarray*}
&\Delta:D_{(k,l)}\hookrightarrow(\Rf\times\mathcal{M}_{k})^l\\
&\left(x_1,\ldots,x_l,[T,P]\right)\mapsto\left((x_1,[T,P]),(x_2,[T,P]),\ldots,(x_l,[T,P])\right)
\end{eqnarray*}
a sort of diagonal map.
\par
A similar method with the infinite dimensional case yields the construction of
the Donaldson polynomial.
We will go on and show that although this integral does not vanish for the
$k=1$ case where there is a discrepancy manifesting in the linking phenomena,
the analogous integrals do vanish for higher $k$.
%%\input{equiv1}
%%\Endoff{equiv1.tex}
%%\input{equiv2}
%%\input{equiv3}
%\input{equiv4}
\section{Equivariant Characteristic Classes}
\subsection{General Theory}
We wish to find an equivariant representative of the $\mu$-form of a
submanifold in $\Rf$. Since this form was built up from a characteristic
class it is necessary for us to consider the theory of equivariant
characteristic classes.\par
Let $p:V\mapping M$ be an equivariant $G$ bundle over the $G$-manifold $M$.
We may form the new bundle $V_G\mapping M_G$ by
\[
V_G=EG\times_G V \mapping EG\times_GM =M_G
\]
and compare the characteristic classes of $V$ with those of $V_G$.\par
Given an equivariant connection $\grad$ on $V$ we may form a connection $\grad_G$
on $V_G$ by pulling back by the projection $EG\times M\mapping M$ and observing
how it descends to the quotient. Now if $K$ is the structure group of $V$, then it is
also the structure group of $V_G$, so given a $K$-invariant polynomial
$P\in\left(\bigodot\goth{k}^*\right)^K$ we may form the characteristic classes
\begin{eqnarray*}
&c^P=P(F(\grad))\in\Form{M}{\mathrm{even}}\\
&c_G^P=P(F(\grad_G))\in\OG{M}{even}.
\end{eqnarray*}
We should like to see how these are related.
\begin{prop1}[Selby \cite{S} p16]\hfill\\
If $f:M\mapping N$ is $G$-equivariant  between $G$-manifolds inducing the map
\[
f_G:EG\times_G M\mapping EG\times_G N
\]
and $p:V\mapping N$ a $G$-equivariant fibre bundle then
\[
f_G^*V_G=(f^*V)_G
\]
\end{prop1}
The proof is easy and can be found in detail in \cite{S}
Choose a base point $e\in EG$ and let $\iota:M\mapping EG\times M$ be the
inclusion $m\mapsto (e,m)$. So for an equivariant vector bundle
$p:V\mapping M$ we have
\begin{eqnarray*}
V &\iso
&\iota^*(V\times EG)\\
&\iso&\iota^* q^*V_G
\end{eqnarray*}
where $q:EG\times M\mapping M_G$ is the quotient map.We must be careful here;
this is an isomorphism of bundles but is not necessarily equivariant.
If we set $r=q\iota:M\mapping M_G$ then $r^*:H_G^*(M)\mapping H^*(M)$, and
further it can also be shown that
\[
r^*c_G^P=c^P
\]
\subsection{de Rham Theory}\label{equivchar}
We now translate the equivariant theory into de Rham formalism. Recall that we
make the identification
\[
\Form{M_G}{*}\mapping\OG{M}{*}=\left(\Form{M}{*}\tensor\complexes[\goth{g}^*]
\right)^G
\]
we can also make the identification
\[
\Form{M_G;V_G}{*}\mapping\left(\Form{M;V}{*}\tensor\complexes[\goth{g}^*]
\right)^G
\]
and call these the equivariant forms on $M$ with values in $V$.
Following \cite{BGV}, given an equivariant connection $\grad$ on $V\mapping M$,
we can form the de Rham version $\conng$ of the corresponding connection on
$V_G\mapping M_G$ by setting
\[
\left(\conng s\right)(\xi)=\grad (s(\xi))-X_\xi\hook s(\xi).
\]
As mentioned in \cite{BGV} pp210-211, we are motivated by
\[
\dg^2s(\xi)+\mathcal{L}_\xi s(\xi)=0
\]
to define the equivariant curvature
\[
\curv(\grad)s(\xi)=\conng^2s(\xi)+\mathcal{L}_\xi^Vs(\xi)
\]
whence
\begin{equation}\label{equivcurv}
\curv(\grad)s(\xi)=F(\grad)s(\xi)-[\grad,X_\xi\hook]s(\xi)+\mathcal{L}_\xi s(\xi)
\end{equation}

where $\mathcal{L}_\xi^V$ is the Lie derivative on $V$ induced by the action of
the vector field $X_\xi$. Now, given a $K$-invariant polynomial $P$ ($K$ being the structure group of $V$), we may obviously form an equivariant
characteristic class
\[
c^P(\xi)=P(\curv(\grad){\xi}\in\OG{M}{\mathrm{even}}.
\]
How does this relate to the corresponding characteristic class of $V/G$ on
$M/G$?
\par
If $G$ does not act freely on $M$, then consider the manifold
$M^*=M\backslash M_0$ of points with trivial stabiliser. Hence $M^*/G$ is a
manifold.
\begin{lem1}\label{equivchar1}
For a vector bundle $V\mapping M^*$
\[
c_G^P(V)=q^*c^P(V/G)+\dg \beta
\]
where $q:M\mapping M/G$ is the quotient map.
\end{lem1}
\begin{prf}
Given an equivariant connection $\grad$ on $V$, we may proceed
as in {section 5.2.3 of \cite{DK}}
%{\thingyrefi}
 to obtain a connection $\grad'$ on $V/G\mapping M^*/G$.
In turn $q^*\grad'$ determines an equivariant and horizontal connection
on $V\mapping M^*$ and hence an equivariant de Rham connection $\conng'$ on
$V_G\mapping M_G^*$. \par
We therefore have two connections on $V_G\mapping M_G^*$ namely $\conng$ and
$\conng'$. By the usual arguments in the theory of characteristic classes
\[
P(\curv(\grad'))=P(\curv(\grad))+\dg \beta.
\]
But $P(\curv(\grad')=q^*P(F(\grad'))$ and defines an equivariantly closed form
on $M^*$. So the result follows.
\end{prf}
\begin{cor1}
Any equivariant characteristic class of a bundle over a hyperK\"ahler manifold
is associated to the characteristic class of the quotient bundle over the
hyperK\"ahler reduction.
\end{cor1}
\subsection{HyperK\"ahler Integration}
We use the results developed in \cite{MU2}. 
\begin{defn1}
Let  ($M$,$\vec{\omega}=\omega_\I\I+\omega_\J\J+\omega_\K\K$) be a hyperK\"ahler manifold which possesses a %\newline
tri-Hamiltonian action of the compact Lie group $G$ and
$\alpha\in\OG{M}{\bullet}$ be an equivariant form. We shall say
that $\alpha$ is associated to
$\alpha_0\in\Form{\mathcal{M}}{\bullet}$ if
\[
\iota^*\alpha=\pi^*\alpha_0+\dg\beta
\]
where $\iota:\vec{\mu}^{-1}(0)\mapping M$ is inclusion and
$\pi:\vec{\mu}^{-1}(0)\mapping\mathcal{M}$ is the quotient map and
$\beta\in\OG{\vec{\mu}^{-1}(0)}{\bullet}$. In the case that
$\alpha$ and $\alpha_0$ are compactly supported, we shall say that
$\alpha$ is compactly associated to $\alpha_0$ if $\alpha$ is
associated to $\alpha_0$ as above, and the form $\beta$ is also
compactly supported.
\end{defn1}
The reduction of $(M,\vec\omega)$ is a manifold $(M/\!\!/\!\!/\!\!/G,\vec{\omega}_0)$.
The compactness of support is a necessity from the fact communicated to the Author
by Roger Bielawski that there are no compact hyperK\"ahler reductions by
a tri-Hamiltonian group. One of the main results of \cite{MU2} (also found
in \cite{MU}) is
\begin{thm1}\label{HYPINTFORM}
If $\eta$ is compactly associated to $\eta_0$, then
\begin{eqnarray*}
& &\int_\mathcal{M}\exp{i\vec{\omega}_0\wdot\vec{\omega}_0}\left(4i\vec{\omega}_0\wdot\vec{\omega}_0+1\right)\eta_0\\
&=&\left(\frac{1}{6\pi i\sqrt{2}}\right)^k\frac{1}{|W|}\oint_M
\exp{i\vec{\omega}\wdot\vec{\omega}+i|{\vec{\mu}}|^2y}\mathrm{Pr}_{ev}\left(z\mapsto\exp{2i{z}{\vec{\mu}}.\vec{\omega}}w(z)^4\eta(z)\right)(\sqrt{y})\\
\end{eqnarray*}
where
\[
\oint_M\alpha(\xi)=\frac{1}{\mathrm{vol}\,(G)}\lim_{t\mapping\infty}
\int_\g\d\xi\exp{-\frac{|\xi|^2}{4t}}\int_M\alpha(\xi),
\]
and
\[
\vec{\omega}\wdot\vec{\omega}=\omega_\I\wedge\omega_\I+\omega_\J\wedge\omega_\J+\omega_\K\wedge\omega_\K,
\]
the operator $\mathrm{Pr}_{ev}:\OG{M}{\bullet}\mapping\OG{M}{\bullet}$ is
\[
\left(\mathrm{Pr}_{ev}(\beta)\right)(z)=\frac{1}{2}\left(\beta(z)+\beta(-z)\right),
\]
$|W|$ is the number of elements of the Weyl Group associated with the maximal torus in
$G$
and 
\[
w(y)=\prod_{\alpha\in\Delta_+}\alpha(y)
\]
is the polynomial formed by the product of the positive roots of $G$.
\end{thm1}
It is this result that allows us to compute the integrals in which we are primarily interested.
\section{Applications to the ADHM construction}
We wish to apply the results on equivariant integration and
localisation for hyperK\"ahler quotients to the ADHM construction.
Here it is better to pass to the ``complex'' version given by
\[
(T,P)\in\mc{k}=\left(i\u{k}\oplus\complexes^k\right)\tensor_\reals\quat
\]
which has the moment map
\[
\vec{\mu}(T,P)=\Im_\quat(T^*T+PP^*)
\]
where $\Im_\quat$ is the complexification of the operation of
taking the quaternionic imaginary part and $*$ means taking the
quaternionic conjugate of the complex adjoint.\par The reason for
this change of approach is that
\[
\widetilde{\mathcal{M}}=\mr{k}{/\!\!/\!\!/\!\!/}\O{k}\iso\mc{k}{/\!\!/\!\!/\!\!/}\U{k}
\]
and $\U{k}$ is connected and has a simpler Lie algebra structure
than $\O{k}$. We may also describe the maximal torus of $\U{k}$
more simply than $\O{k}$ and we will be using this to assist us in
our localisation. This does not affect any of our previous
results.\par
\subsection{The Equivariant Euler Class}
Our first priority is to work out the fixed point set of the
action of the maximal torus
\[
\mathbb{T}^k\subset\U{k}.
\]
To do this, we take the decomposition
\[
\mathbb{T}^k=\mathbb{T}_1\times\ldots\times\mathbb{T}_k
\]
where
\[
\mathbb{T}_j=\left\{\left. \left(
\begin{array}{ccc}
 \id & 0 & 0\\
   0 & e^{i\theta_j} & 0\\
   0 & 0 & \id
\end{array}
\right) \right\vert \theta_j\in\reals \right\}.
\]
By finding the fixed set of $\mathbb{T}_i$ and the equivariant
Euler class of its normal bundle in $\mc{k}$ we will be able to
apply our inductive formula.\par Let $\xi_j\in\mathbb{T}_i$ be the
generator
\[
\xi_j= \left(
\begin{array}{ccc}
 0& 0 & 0 \\
 0& i & 0 \\
 0& 0 & 0
\end{array}
\right).
\]
Then the vector field $X_{\xi_j}$ on $\mc{k}$ is given by
\begin{eqnarray*}
X_{\xi_j}(T,P)&=&\left([\xi_j,T],\xi_j P\right)\\
              &=&
\left( i \left[
\begin{array}{ccc}
 0 & -T_1 & 0\\
T_1^* & 0    & T_4 \\
0  & -T_4^* & 0
\end{array}
\right], \left[
\begin{array}{c}
0 \\
\vdots\\
iP_j\\
0\\
\vdots\\
0
\end{array}
\right] \right),
\end{eqnarray*}
where
\[
T=\left(
\begin{array}{ccc}
 T' & T_1 & T_2 \\
 T_1^*& T'' & T_4 \\
 T_2^*& T_4^* & T'''
\end{array}
\right), \hbox{ and } P= \left(
\begin{array}{c}
P_1 \\
\vdots\\
P_k
\end{array}
\right).
\]
Let $\mathrm{Sing}_j^k$ consist of the points where $X_{\xi_j}$
vanishes, i.e
\[
\mathrm{Sing}_j^k:=\left\{\left(\left[\begin{array}{ccc}
 T' & 0 & T_2 \\
 0& T'' &0 \\
 T_2^*& 0 & T'''
\end{array}
\right], \left[
\begin{array}{c}
P' \\
0\\
P''
\end{array}
\right] \right)\in \mc{k} \right\}.
\]
\begin{thm1}\label{euler}
The $\sphere{1}$-equivariant Euler class $e$ of the normal bundle
of $\mathrm{Sing}_j^k$ in $\mc{k}$ is given by
\[
e(\lambda)=\left(\frac{\lambda}{2\pi}\right)^{4k}.
\]
\end{thm1}
\begin{prf}
Set
$\mathcal{V}_j=\mathcal{V}(\mathrm{Sing}_j^k\hookrightarrow\mc{k})$,
the normal bundle of $\mathrm{Sing}_j^k$ in $\mc{k}$. We notice
that $\mathrm{Sing}_j^k$ is a vector subspace of $\mc{k}$, so
\begin{eqnarray*}
\tang{\mathrm{Sing}_j^k}{(T,P)}&\iso&\mathrm{Sing}_j^k,\\
\tang{\mc{k}}{(T,P)}&\iso&\mc{k},\\
\end{eqnarray*}
canonically. So
\[
\left(\mathcal{V}_j\right)_{(T,P)}=\left(\tang{\mathrm{Sing}_j^k}{(T,P)}\right)^\perp=(\mathrm{Sing}_j^k)^\perp.
\]
It can be shown that
\begin{eqnarray*}
(\mathrm{Sing}_j^k)^\perp&=&\{ X_{\xi_j}(T,P)\vert (T,P)\in\mc{k}\}\\
             &=&\left\{\left.\left(\left[
\begin{array}{ccc}
 0 & T_1 & 0\\
-T_1^* & 0    & T_4 \\
0  & -T_4^* & 0
\end{array}
\right], \left[
\begin{array}{c}
0 \\
\vdots\\
P_j\\
0\\
\vdots\\
0
\end{array}
\right]\right)\right\vert \begin{array}{rcl}
T_1&\in&\complexes^{j-1}\tensor_\reals\quat,\\
T_4&\in&\complexes^{k-j}\tensor_\reals\quat,\\
P_j&\in&\complexes\tensor_\reals\quat
\end{array}
\right\}.
\end{eqnarray*}
Since $\mathcal{V}_j$ is an equivariant bundle, and the de-Rham
operator is a perfectly good equivariant connection on
$\mathcal{V}_j$, we see automatically from (\ref{equivcurv}) that
the equivariant curvature
\[
F_{\goth{g}}(\d)(\lambda\xi_j)=\mathcal{L}_{\lambda\xi_j}^{\mathcal{V}_j}.
\]
Now,
\[
\mathcal{L}_{\lambda\xi_j}^{\mathcal{V}_j}\left(\left[
\begin{array}{ccc}
 0 & T_1 & 0\\
-T_1^* & 0    & T_4 \\
0  & -T_4^* & 0
\end{array}
\right], \left[
\begin{array}{c}
0 \\
\vdots\\
P_j\\
0\\
\vdots\\
0
\end{array}
\right]\right) =\lambda i \left(\left[
\begin{array}{ccc}
 0 &  -T_1 & 0\\
T_1^* & 0    & -T_4 \\
0  & T_4^* & 0
\end{array}
\right], \left[
\begin{array}{c}
0 \\
\vdots\\
P_j\\
0\\
\vdots\\
0
\end{array}
\right]\right)
\]
More simply
\[
\mathcal{L}_{\lambda\xi_j}^{\mathcal{V}_j}(T_1,T_4,P_j)=(-i\lambda
T_1,-i\lambda T_4,i\lambda P_j).
\]
From this we are able to deduce that
\begin{eqnarray*}
e_j(\lambda)&=&\mathrm{Pfaff}(\frac{1}{2\pi}F_{\goth{g}}(\lambda\xi_j))\\
            &=&\left(\frac{-\lambda}{2\pi}\right)^{4(k-1)}\left(\frac{\lambda}{2\pi}\right)^4\\
            &=&\left(\frac{\lambda}{2\pi}\right)^{4k}.
\end{eqnarray*}
%Now, for $(T,P)\in \mathrm{Sing}_j^k$, we have
%\[
%\begin{array}{ll}
%& T^*T+PP^*=\\
%& \\
%&\left(
%\begin{array}{ccc}
%(T')^\dag (T')+ (T_2)^\dag(T_2^*)+(P')(P')^* & 0 & (T')^\dag(T_2)+(T_2)^\dag (T''')+(P')(P'')^*\\
%0 & |T''|^2 &0\\
%(T_2)^\dag(T')+(T''')^\dag (T_2^*)+(P'')(P')^* & 0& (T''')^\dag (T''')+
%(T_2)^\dag(T_2)+(P'')(P'')^*
%\end{array}
%\right)
%\end{array}
%\]
%where $\dag$ denotes quaternionic conjugation.
%Now, if we choose
%\[
%\vec{\zeta}=\left(
%\begin{array}{ccc}
%0 & 0 & 0\\
%0 & a & 0\\
%0 & 0 & 0
%\end{array}
%\right)
%\]
%for $a\in\sP{1}$ non-zero,
%we observe that $\inv{\vec{\mu}}(\zeta)$ cannot intersect $\mathrm{Sing}_j^k$.
\end{prf}
\subsection{The Universal 2nd Chern Class}
In order to find the Universal equivariant 2nd Chern class, we
apply the theory in \ref{equivchar} to the connection $\hcon=v^*\circ\hd\circ
v$ on the the bundle $\hE\mapping\Rf\times\mc{k}$ 
%that we created in {\thingyrefii}. Recall that
which has curvature
\[
F(\hcon)=v^*\hd\R^*\wedge\F\hd\R v.
\]
Using this it can be quite easily shown that for some section
local section
$s\in\sect{\mc{k}}{\hE}$ and $\xi\in\u{k}%,\vec{q}\in\sP{1}
$, we have
\begin{eqnarray*}
\hcon_{X_{\xi%\vec{q}
}}s&=&X_{\xi%\vec{q}
}\hook v^*\hd vs\\
                         &=&v^*\hd(vs)(X_{\xi%\vec{q}
})\\
                         &=&v^*\Diff{t}\restrict{t=0}\left(\phi_t^{(\xi%\vec{q}
)}\right)^*vs\\
                         &=&v^*\left(\Diff{t}\restrict{t=0}\left(\phi_t^{(\xi%\vec{q}
)}\right)^*v\right)s
+\Diff{t}\restrict{t=0}\left(\phi_t^{\xi}\right)^*s\\
&=&v^*\left(\begin{array}{cc}
\xi & 0\\
0   & 0%-\vec{q}

\end{array}\right)vs
+\mathcal{L}_{\xi}s.
\end{eqnarray*}
Hence
\begin{eqnarray*}
\curv(\hcon)(\xi%\vec{q}
)&=&F(\hcon)-v^*\left(\begin{array}{cc}
\xi & 0\\
0   & 0%-\vec{q}

\end{array}\right)v\\
&=&v^*\left(\hd\R^*\wedge\F\hd\R-\left(\begin{array}{cc}
\xi & 0\\
0   & 0%-\vec{q}

\end{array}\right)\right)
v
\end{eqnarray*}
and obviously our representative of the equivariant universal
Chern class will be
\[
c^k(\xi%\vec{q}
)=\frac{1}{4\pi^2}\mathrm{tr}\left(\curv(\hcon)(\xi%\vec{q}
)\wedge\curv(\hcon)(\xi%\vec{q}
)\right)\in\Omega_{\mathbb{T}^k%\SP{1}
}^4(\Rf\times\mc{k}).
\]
However, $c^k$ is integrable but not smooth on the space of
reducibles. To understand $c^k$ on the space of reducibles, if we
look at $c^k_{(x,T,P)}$ for
\[
T=\left(\begin{array}{cc}
T_0 & 0\\
0& T_1
\end{array}\right)
,\quad P=\left(\begin{array}{c}
0\\
P_1
\end{array}\right)
\]
where $T_0\in\quat$, $(T_1,P_1)\in\mc{k-1}$ and $x\neq T_0$, then
it is straightforward to show that
\[
c^k_{(x,T,P)}(\xi%\vec{q}
)=c^{k-1}_{(x,T_1,P_1)}(\xi_1%,\vec{q}
)
\]
where $\xi_1$ is the element of $\mathrm{Lie}\mathbb{T}^{k-1}$ got
from $\xi$ by removing the first row and column. But,
\[
\int_{\Rf}c^k_{(x,t,p)}(\xi%\vec{q}
)\d x=k
\]
for all irreducible $(t,p)\in\mc{k}$.
\begin{lem1}\label{deltachern}
For the above $(T,P)$
\[
\iota_{1}^*c^k_{(x,T,P)}(\xi%\vec{q}
)=c^{k-1}_{(x,T_1,P_1)}(\xi_1%,\vec{q}
)+\delta(x-T_0)\d(x-T_0)_1\wedge\d(x-T_0)_2\wedge\d(x-T_0)_3\wedge\d(x-T_0)_4.
\]
where $\iota_1:\mathrm{Sing}_1^k\hookrightarrow\mc{k}$ is
inclusion and if
\[
\xi=\left(
\begin{array}{cccc}
\xi_{11} &0&       &  \\
 0     &\xi_{22} & &  \\
       & &\ddots &  \\
       & &       &\xi_{kk}
\end{array}
\right)
\]
we let
\[
\xi_1=\left(
\begin{array}{ccc}
\xi_{22} & &  \\
  &\ddots &  \\
  &       &\xi_{kk}
\end{array}
\right).
\]
\end{lem1}
\begin{prf}
Set
\[
L(\xi%\vec{q}
)=\left(\begin{array}{cc}
\xi & 0\\
0   & 0%-\vec{q}

\end{array}\right).
\]
The splitting of the matrix $L(\xi%\vec{q}
)$ is in terms of the splitting of
$\complexes^k\tensor_\reals\quat\oplus\complexes^2$. We have
another splitting to consider here due to separating $T$ into
$T_0$ and $T_1$. Let
\[
L(\xi%\vec{q}
)=\left(\begin{array}{cc}
\xi_{11} & 0\\
0   & L_1(\xi_1%,\vec{q}
)
\end{array}\right),
\]
where $L_1(\xi_1%,\vec{q}
)$ plays same r\^ole as $L$ for $\mc{k-1}$.\par We have
\[
\R_{(x,T,P)}=\left(
\begin{array}{ccc}
(T_0-x)^* &0 &0\\
0         &(T_1-x\id)^*&P_1
\end{array}
\right) = \left(
\begin{array}{cc}
-\tilde{x}_0^* & 0\\
 0          & \R_1
\end{array}
\right)
\]
where $\tilde{x_0}=x-T_0$ and $\R_1=\R_{(x,T_1,P_1)}$. Now
\[
\F=(\R\R^*)^{-1}=\left(
\begin{array}{cc}
|\tilde{x}_0|^{-2} & 0\\
0 & \F_1
\end{array}
\right)
\]
where $\F_1=(\R_1\R_1^*)^{-1}$. Since there is a singularity when $x=T_0$, we make a small adjustment depending on a
parameter $\rho$ which we will shrink to 0.\par Set
\[
\F_\rho=\left(
\begin{array}{cc}
\frac{1}{(\rho^2+|\tilde{x}_0|^{2})} & 0\\
0 & \F_1
\end{array}
\right).
\]
Define
\[
\varpi_\rho=\id-\R^*\F_\rho\R= \left(
\begin{array}{cc}
\frac{\rho^2}{(\rho^2+|\tilde{x}_0|^{2})} & 0\\
0 & \varpi_1
\end{array}
\right)
\]
where $\varpi_1=\id-\R_1^*\F_1\R_1$. We note that
\[
\iota_1^*\d\R=\left(
\begin{array}{cc}
-\hd\tilde{x}_0^* & 0\\
 0          & \hd\R_1
\end{array}
\right).
\]
So we have
\begin{eqnarray*}
& &\iota_1^*c^{k}_{(x,T,P)}(\xi)\\
%\vec{q})\\
 &=&\lim_{\rho\mapping 0}\frac{1}{4\pi^2}\Re\mathrm{tr}
\left[\left((\iota_1^*\hd\R^*\wedge\F_\rho\iota_1^*\hd\R\varpi-L(\xi%\vec{q}
)\varpi\right)^2\right] \\
%&=&\lim_{\rho\mapping 0}\frac{1}{4\pi^2}\Re\mathrm{tr} \left[
%\left(
%\begin{array}{cc}
%\frac{\rho^2}{(\rho^2+|\tilde{x}_0|^2)^2}\hd{\tilde{x}_0}\wedge\hd{\tilde{x}_0}^*-\frac{\rho^2}{(\rho^2+|\tilde{x}_0|^2)}\xi_{11}& 0 \\
%0& (\hd\R_1^*\wedge\F_1\hd\R_1-L_1(\xi_1%,\vec{q}
%))\varpi_1
%\end{array}
%\right)^2
%\right]\\
&=& \lim_{\rho\mapping
0}\frac{1}{4\pi^2}\Re\left(\frac{\rho^4}{(\rho^2+|\tilde{x}_0|^2)^4}\hd{\tilde{x}_0}\wedge\hd{\tilde{x}_0}^*\wedge\hd{\tilde{x}_0}\wedge\hd{\tilde{x}_0}^*\right.\\
&-&\left.\frac{2\rho^4}{(\rho^2+|\tilde{x}_0|^2)^3}\hd{\tilde{x}_0}\wedge\hd{\tilde{x}_0}^*\xi_{11}+\frac{\rho^4}{(\rho^2+|\tilde{x}_0|^2)}\xi_{11}^2
\right)\\
&+&c_{(x,T_1,P_1)}^{k-1}(\xi_1%,\vec{q}
)\\
&=&\lim_{\rho\mapping 0}\frac{1}{4\pi^2}\Re\left(\frac{24\rho^4}{(\rho^2+|\tilde{x}_0|^2)^4}\hd(\tilde{x}_0)_1\wedge\hd(\tilde{x}_0)_2\wedge\hd(\tilde{x}_0)_3\wedge\hd(\tilde{x}_0)_4\right.\\
&-&\left.\frac{2\rho^4}{(\rho^2+|\tilde{x}_0|^2)^3}\hd{\tilde{x}_0}\wedge\hd{\tilde{x}_0}^*\xi_{11}+\frac{\rho^4}{(\rho^2+|\tilde{x}_0|^2)}\xi_{11}^2
\right)\\
&+&c_{(x,T_1,P_1)}^{k-1}(\xi_1%,\vec{q}
).\\
\end{eqnarray*}
Now the terms of the form
\[
\frac{\rho^4}{\left(\rho^2+|\tilde{x}_0|^2\right)^n}
\]
are distributional in the limit as $\rho_{ii}\mapping0$, so it is
worth integrating them against a compactly supported test function
$f:\Rf\mapping\reals$ that is, calculating
\[
\lim_{\rho\mapping
0}\int_{\Rf}\frac{\rho^4}{\left(\rho^2+|y|^2\right)^n}f(y)\d y.
\]
It isn't hard to show that as distributions
\begin{eqnarray*}
&\displaystyle{\lim_{\rho\mapping
0}\frac{\rho^4}{\left(\rho^2+|y|^2\right)^4}=\frac{\pi^2}{6}\delta(y)},\\
&\displaystyle{\lim_{\rho\mapping
0}\frac{\rho^4}{\left(\rho^2+|y|^2\right)^3}}=0
\end{eqnarray*}
and after a little work is is possible to show that
\[
\lim_{\rho\mapping 0}\frac{\rho^4}{\left(\rho^2+|y|^2\right)^2}=0
\]
as a distribution. We are therefore left with the result that
\[
\iota_{1}^*c^k_{(x,T,P)}(\xi%\vec{q}
)=c^{k-1}_{(x,T_1,P_1)}(\xi_1%,\vec{q}
)+\delta(x-T_0)\d(x-T_0)_1\wedge\d(x-T_0)_2\wedge\d(x-T_0)_3\wedge\d(x-T_0)_4.
\]
as stated.
\end{prf}
For a submanifold $\Sigma\subset\Rf$, define
\[
\mu^k(\Sigma)_{(T,P)}(\xi%\vec{q}
)=\int_{x\in\Sigma} c^k_{(x,T,P)}(\xi)=\int_{x\in\Rf}
c^k_{(x,T,P)}(\xi)\wedge PD(\Sigma)_x.
\]
\begin{cor1}\label{redmu} With $T,P,\iota_1$ as above, we have
\[
\iota_1^*\mu^k(\Sigma)_{(x,T,P)}(\xi%\vec{q}
)=\mu^{k-l}(\Sigma)_{(T_1,P_1)}(\xi_1%,\vec{q}
)+\mathrm{PD}(\Sigma)_{T_0}.
\]
\end{cor1}
\begin{cor1}\label{bdy}
If $\iota$ is the inclusion of the set of points of $\mc{k}$ fixed
under $\mathbb{T}^k$, and $p_i:\mc{k}\mapping\Rf$ is the
projection $T\mapsto T_{ii}$ then
\[
\iota^*\mu^k(\Sigma)(\xi%\vec{q}
)=\sum_{i=1}^k p_i^*\mathrm{PD}(\Sigma).
\]
\end{cor1}
\begin{cor1}\label{fixedmu}
If
$\mu^k(\Sigma)\in\Form{\frac{\mc{k}}{\U{k}\SP{1}}}{4-\dim\sigma}$
then for
%\[
%T=\left(\begin{array}{cc} T_0 & 0\\
%0 & T_1 \end{array}\right), \]
%\[ P= \left(
%\begin{array}{c}
%0\\
%P_1\\
%\end{array}
%\right),
%\]
$(T,P)$ and $\iota_1$ as above we have
\[
\iota_1^*\pi^*p^*\mu^k(\Sigma)_{(T,P)}=\pi^*p^*\mu^{k-1}(\Sigma)_{(T_1,P_1)}
+p_1^*\mathrm{PD}(\Sigma)_{T_0}+\pi^*p^*\hd\gamma(\Sigma).
\]
for some
$\gamma(\Sigma)\in\Form{\frac{\mc{k}}{\U{k}\SP{1}}}{3-\dim\Sigma}$,
and $p_1$ the projection as in Corollary \ref{bdy}.
\end{cor1}\par
\begin{prf}
Now, we have a chain of quotients
\[
\begin{array}{ccccc}
\mc{k}&\stackrel{\pi}{\mapping}&\frac{\mc{k}}{\U{k}}&\stackrel{p}{\mapping}&\frac{\mc{k}}{\U{k}\SP{1}}\\
      &                        &\cup                &                      & \cup\\
      &                        & \widetilde{\mathcal{M}}_k&                & \mathcal{M}_k\\
\end{array}.
\]
So the form $\pi^*p^*\mu^k(\Sigma)$ is now a closed, $\U{k}$-basic
$4-\dim\Sigma$ degree form on $\mc{k}$. Since
$\mu^k(\Sigma)(\xi)\in\Omega_{\U{k}}(\mc{k})$ was formed using
equivariant Chern-Weil theory and by Lemma \ref{equivchar1}, we
have
\[
\mu^k(\Sigma)(\xi)=\pi^*p^*\mu^k(\Sigma)+\int_\Sigma\dg\beta^k(\xi).
\]
Also by Corollary \ref{redmu} we have
\[
\iota_1^*\mu^k(\Sigma)_{(x,T,P)}(\xi)=\mu^{k-l}(\Sigma)_{(T_1,P_1)}(\xi_1)+\mathrm{PD}(\Sigma)_{T_0}.
\]
Thus
\[
\iota_1^*\pi^*p^*\mu^k(\Sigma)=\pi^*p^*\mu^{k-1}(\Sigma)_{(T_1,P_1)}
+\mathrm{PD}(\Sigma)_{T_0}+\int_\Sigma\dg\left(\beta^{k-1}-\iota_1^*\beta^k\right)(\xi).
\]
The left hand side is independent of $\xi$, so
$\dg\left(\beta^{k-1}-\iota_1^*\beta^k\right)(\xi)$ is an exact,
$\U{k}$- basic form, thus the right hand side is in the same de
Rham cohomology class as the left. The result follows since
$\pi^*p^*\mu^k(\Sigma)$, $\pi^*p^*\mu^{k-1}(\Sigma)$ and the
Poincar\'e dual are $\U{k}\SP{1}$-basic. By construction,
$\iota_1^*\dg\beta^k$ does not depend on $T_0$ and thus agrees
with $\dg\beta^{k-1}$ giving the result.
%differ by a form in the
%variable $T_0$, thus giving a change of cohomology class for the
%Poincar\'e dual of $\Sigma$. Setting
%$$\gamma(\Sigma)=\mathrm{PD}(\Sigma)
%+\int_\Sigma(\dg\left(\beta^{k-1}-\iota_1^*\beta^k\right)(\xi)$$
%gives the result.
\end{prf}
\subsection{Integrability of the Donaldson $\mu$ map}

What is not altogether clear is that the form representing the
Donaldson polynomial is actually integrable. Indeed there are
various technicalities in forming the these polynomials that
Donaldson and Kronheimer discuss in Chapter 9 of \cite{DK}. Our
approach will be from a functional analytic viewpoint.
\begin{defn1}
Given pairwise disjoint, compact submanifolds $\Sigma_1,\ldots\,\Sigma_l$ of $\Rf$ we
define the Donaldson functional on compactly supported functions
of $\frac{\mc{k}}{\U{k}\SP{1}}$ by
\[
\mathrm{Don}_k(\Sigma_1,\ldots,\Sigma_l)(\phi)
=\int_{\mathcal{M}_k}\phi\,\mu(\Sigma_1)\wedge\ldots\wedge\mu(\Sigma_l).
\]
\end{defn1}
This is certainly well defined; the above argument shows that the
representatives of the $\mu$ classes on the reducible space are
distributional in nature and integrable, and thus the integral
exists for any compactly supported function $\phi$.

%We will however make a particular choice for $\phi$. Let $K$ be a
%compact subset of $\Rf$ which contains $\bigcup_{i=1}^l\Sigma_i$.
%We require $\phi$ to be identically 1 on the space
%\[
%\frac{K\times\ldots\times K}{\U{k}\SP{1}}
%\]
%regarded as a subset of the quotient of the fixed set of the
%action.
%\begin{defn1}
%We call such a $\phi$ an amicable test function.
%\end{defn}
%\input{equiv6a.tex}
%\subsection{Dealing with the $\SO{3}$-action}
\subsection{Computing the integrals}
We must, however, make a slight alteration to the situation since
the action of $\U{k}$ is not free on $\vec{\mu}^{-1}(0)$. Instead
we choose $\vec{\zeta}_0\in\Im\quat$ and take the moment map to be
\[
\vec{\mu}_{\zeta_0}(T,P)=\Im(T^*T+PP^*)-\vec{\zeta_0}\id.
\]
%\subsection{Dealing with the $\SO{3}$ action}
We have to decide on how best to approach the integration.

Let $\Sigma_1,\ldots,\Sigma_l$ be pairwise disjoint, compact submanifolds of $\Rf$ of
dimensions $d_1,\ldots,d_l$ respectively such that
\[
\sum_{i=1}^{l}(4-d_i)=8k-3,
\]
that is
\[
\sum_{i=1}^{l}d_i=4l-8k+3.
\]
Then $\alpha =\mu(\Sigma_1)\wedge\ldots\wedge\mu(\Sigma_l) $  is
represented by
 a  form of top degree on ${\mathcal{M}_k}$.
%and integrable on all of $\mathcal{M}_k$.

Now the de-Rham operator $\hd$ splits on $\Rf\times\mathcal{M}_k$
\[
\hd=\d+\delta
\]
where $\delta$ is the de Rham differential on $\mathcal{M}$. Since
both $\U{k}$ and $\SP{1}$ act trivially on $\Rf$, we see that
\[
\dg\eta(\gamma)=\d\eta(\gamma)+\delta\eta(\gamma)-X_\gamma\hook\eta(\gamma)
=\hd\eta(\gamma)-X_\gamma\hook\eta(\gamma),
\]
for any $\eta\in\OG{\Rf\times\mc{k}}{\bullet}$, where following
the notation in Chapter 2, we reserve $\d$ for the de-Rham
differential on $\Rf$ and $\delta$ the differential on $\mc{k}$
and the total differential $\hd=\d+\delta$. Now suppose without
loss of generality that $\Sigma_1$ is not a point, and that
$\Xi_1$ is a Seifert surface spanning $\Sigma_1$, i.e
$\partial\Xi_1=\Sigma_1$. Then
\begin{eqnarray*}
0&=&\int_{\Xi_1} \hd c^k\\
 &=&\int_{\Xi_1} \d c^k +\int_{\Xi_1}\delta c^k\\
 &=&\int_{\Sigma_1} c^k+\delta\int_{\Xi_1}c^k
\end{eqnarray*}
i.e
\begin{equation}\label{exactness}
\int_{\Sigma_1} c^k=-\delta\int_{\Xi_1}c^k.
\end{equation}

Hence we may take
\[
\beta=-\int_{\Xi_1}c^k\wedge\int_{\Sigma_2}c^k\wedge
\ldots\wedge\int_{\Sigma_l}c^k
\]
and
\[
\alpha=\delta\beta.
\]
%Let us set
%\[
%\mathrm{Don}_k(\Sigma_1,\ldots,\Sigma_l)=\int_{\mathcal{M}_k}\mu(\Sigma_1)\wedge\ldots\wedge\mu(\Sigma_l)
%\]
So as we saw above,
\begin{eqnarray*}
& &\mathrm{Don}_k(\Sigma_1,\ldots,\Sigma_l)(\phi)\\
&=&\int_{\mathcal{M}_k}\phi\,\mu(\Sigma_1)\wedge\ldots\wedge\mu(\Sigma_l)\\
&=&\int_{\widetilde{\mathcal{M}}_k}p^*\phi\,
p^*\hd\left(\mu(\Xi_1)\wedge\ldots\wedge\mu(\Sigma_l)\right)\wedge\Theta\\
&=&\int_{\widetilde{\mathcal{M}}_k}
\exp{i\vec{\omega_0}\wdot\vec{\omega_0}}\left(4i\vec{\omega_0}\wdot\vec{\omega_0}+1\right)p^*\phi\,p^*\d\left(\mu(\Xi_1)\wedge\ldots\wedge\mu(\Sigma_l)\right)\wedge\Theta\\
\end{eqnarray*}
since $$p^*\phi
p^*\d\left(\mu(\Xi_1)\wedge\ldots\wedge\mu(\Sigma_l)\right)\wedge\Theta$$
already has maximal degree. This form is associated to
$$\eta=\pi^*\left(p^*\phi p^*\hd\left(\mu(\Xi_1)\wedge\ldots\wedge\mu(\Sigma_l)\right)\wedge\Theta\right)
\in\Omega_{\U{k}}^{8k}\left(\mc{k}\right),$$ which is basic and de
Rham closed by construction and compactly supported.\par We can
now use Theorem \ref{HYPINTFORM} applied to this form.
\begin{eqnarray*}
& &\mathrm{Don}_k(\Sigma_1,\ldots,\Sigma_l)(\phi)\\
&=&\left(\frac{1}{6\pi
i\sqrt{2}}\right)^k\frac{1}{|S_k|}\oint_{\mc{k}}
\exp{i\vec{\omega}\wdot\vec{\omega}+i|{\vec{\mu}_{\zeta_0}}|^2y}\mathrm{Pr}_{ev}\left(z\mapsto\exp{2i{z}{\vec{\mu}_{\zeta_0}}.\vec{\omega}}w(z)^4\eta\right)(\sqrt{y}),\\
&=&\left(\frac{1}{6\pi
i\sqrt{2}}\right)^k\frac{1}{k!}\oint_{\mc{k}}
\exp{i\vec{\omega}\wdot\vec{\omega}+i|{\vec{\mu}_{\zeta_0}}|^2y}\mathrm{Pr}_{ev}\left(z\mapsto\exp{2i{z}{\vec{\mu}_{\zeta_0}}.\vec{\omega}}w(z)^4\right)(\sqrt{y})\eta,\\
\end{eqnarray*}
and use the localisation theorem to prove
\begin{thm1}\label{biggie}
\[
\mathrm{Don}_k(\Sigma_1,\ldots,\Sigma_l)(\phi)=\lambda(|\vec{\zeta}_0|)\P(\Sigma_1,\ldots,\Sigma_l)(\phi)
\]
for $\lambda$ a suitable polynomial, and $\P$ a topological quantity
depending on the arrangements of the $\Sigma_i$ in $\Rf$ and upon
the test function $\phi$ .
\end{thm1}
\begin{prf}
We localise the integral with respect to the $(k-1)$-torus in
stead of the $k$-torus since there is a problem with the form
$\Theta$ at $|P|=0$ which is the fixed set of the full $k$-torus .
We hope to be able to express the integral then in terms of the
Donaldson polynomials for charge $k=1$. We are not interested in
the constant multiples that occur here, so they will be largely
forgotten.
\begin{eqnarray*}
& &\mathrm{Don}_k(\Sigma_1,\ldots,\Sigma_l)(\phi)\\
&=&\left(\frac{1}{6\pi
i\sqrt{2}}\right)^k\frac{1}{k!}\oint_{\mc{k}}
\exp{i\vec{\omega}\wdot\vec{\omega}+i|{\vec{\mu}_{\zeta_0}}|^2y}\mathrm{Pr}_{ev}\left(z\mapsto\exp{2i{z}{\vec{\mu}_{\zeta_0}}.\vec{\omega}}w(z)^4\right)(\sqrt{y})\eta,\\
&=&const\oint_{(\Rf)^{k-1}\times\mc{1}}\mathrm{Coeff}_{y_1^{-2}\ldots
y_{k-1}^{-2}}\left[\frac{
\iota^*\exp{i\vec{\omega}\wdot\vec{\omega}+i\sum_{\nu=1}^k\left(|\vec{\zeta}_0|^2y_\nu^2+2{y_\nu}{\vec{\zeta}_0}.\vec{\omega}\right)}w(y)^4
}{y_1^{4k}y_2^{4k-4}\ldots y_{k-1}^8}
 \right]\iota^*\eta.
\end{eqnarray*}
Now restricted to the fixed set of the $(k-1)$-torus we have by
Corollary \ref{fixedmu}
\begin{eqnarray*}
%& &
\iota^*\eta
%\\
%&=&\iota^*\pi^*\left(p^*\phi
%p^*\hd\left(\mu(\Xi_1)\wedge\ldots\wedge\mu(\Sigma_l)\right)\wedge\Theta\right)\\
%&=&\hd\left(\left(pi^*p^*\mu^{1}(\Xi_1)
%+\sum_{m=1}^{k-1}p_i^*PD(\Xi_1) +\hdp^*\gamma
%\right)\right.\\
%&\wedge&\left.\prod_{j=2}^l\left(\pi^*p^*\mu^{1}(\Sigma_j)
%+\sum_{m=1}^{k-1}p_i^*PD(\Sigma_j)+\hd
%p^*\gamma(\Sigma)\right)\right)\wedge\iota^*\pi^*\Theta\\
%&=&\hd\left(\left(pi^*p^*\mu^{1}(\Xi_1)
%+\sum_{m=1}^{k-1}p_i^*PD(\Xi_1)\right)\right.\\
%&\wedge&\left.\prod_{j=2}^l\left(\pi^*p^*\mu^{1}(\Sigma_j)
%+\sum_{m=1}^{k-1}p_i^*PD(\Sigma_j)\right)\right)\wedge\iota^*\pi^*\Theta\\
&=&\iota^*\prod_{j=1}^l\left(\pi^*p^*\mu^{1}(\Sigma_j)
+\sum_{m=1}^{k-1}p_i^*PD(\Sigma_j)
%+\hd\gamma(\Sigma_j)
\right)\wedge\iota^*\pi^*\Theta\\
\end{eqnarray*}
where $p_i$ are the projections described in Corollary \ref{bdy}.
Thus the form of the Donaldson polynomial can be seen certainly as
the  product of a polynomial in $|\vec{\zeta}_0|$ and a sum of
integrals of the form
\begin{eqnarray*}
& &
\int_{\mc{1}\times(\Rf)^{k-1}}\left(\iota^*\pi^*p^*\phi\right)\,\pi^*\left(\Theta\wedge\prod_{i\in
I}p^*\mu^1(\Sigma_i)\right)\wedge\prod_{B}\prod_{b\in
B}\mathrm{PD}(\Sigma_b)\\
&=&\int_{\mathcal{M}_1\times(\Rf)^{k-1}}\left(\iota^*\phi\right)\,\left(\prod_{i\in
I}\mu^1(\Sigma_i))\right)\wedge\prod_{B}\prod_{b\in
B}\mathrm{PD}(\Sigma_b)\\
\end{eqnarray*}
where $I$ is a subset of $\{1,\ldots,l\}$, and $B$ runs over all
subsets that form a partition of $\{1,\ldots,l\}\backslash I$.
\end{prf}
%\input{equiv7}
%\section{Linking in Donaldson invariants}
%Let us see how the results of the last section yield the linking theory.
\subsection{$k=1$ Revisited from the Topological Viewpoint}
We present a slightly different approach to the theory of $k=1$. Here we use
the Poincar\'e duality property detailed in Donaldson's paper \cite{DO}. We
state it in the version it appears in \cite{DK}.
\begin{lem1}[Corollary 5.3.3 of \cite{DK} p199]
Let $X$ be a simply connected Riemannian 4-manifold and $E\mapping X$ have $c_2(E)=1$, and let $\tau:X\mapping
\mathcal{B}^*_{X,E}$ be any map into the space of the gauge equivalence
classes of irreducible connections on $E$ with the property that for all $x$, the
connection $\tau(x)$ is flat and trivial outside some ball of finite diameter
centred on $x$. Then the composite
\[
\mathrm{H}_2(X;\integers)\stackrel{\mu}{\mapping}\mathrm{H}^2(\mathcal{B}^*_{X,E})\stackrel{{\tau^*}}{\mapping}\mathrm{H}^2(X;\integers)
\]
is the Poincar\'e duality isomorphism.
\end{lem1}
\begin{defn1}
We will call such a $\tau$ a tractator
\end{defn1}
This can be proved at the level of forms to show that in the case of $X=\sphere{4}$ we have the following
\begin{lem1}\label{Scrunch}
For a tractator $\tau:\sphere{4}\mapping\mathcal{B}_{\sphere{4},E}^*$
, we have for each submanifold $\Sigma$
of $\sphere{4}$
\[
\int_\Sigma\iota_\Sigma^*\alpha=\int_{\sphere{4}}\alpha\wedge\tau^*\mu(\Sigma)
\]
for any $\alpha\in\Form{\sphere{4}}{\dim\Sigma}$, where $\iota_\Sigma$ is the
inclusion of $\Sigma$ in $\sphere{4}$.
\end{lem1}
Now, for $\eps>0$ let
\[
{\mc{1}}_\eps=\{\left.(T,P)\in\mc{1}\right| |P|\ge\eps\}
\]
and
\[
\mathcal{M}_\eps=({\mc{1}}_\eps/\!\!/\!\!/\!\!/\U{1})/\SP{1}.
\]
This is a manifold with boundary.\par
Let $\Sigma_1,\ldots,\Sigma_l$ be pairwise disjoint, compact submanifolds of $\sphere{4}$ with dimensions
$d_1,\ldots,d_l$ respectively such that
\[
d_1+\ldots+d_l=4l-5.
\]
Suppose w.l.o.g that $\Sigma_1$ is not a point and let $\Xi_1$ be a Seifert
manifold for it. then from above we know that
\[
\mu(\Sigma_1)\wedge\ldots\wedge\mu(\Sigma_l)=\d\left(\mu(\Xi_1)\wedge\mu(\Sigma_2)\wedge\ldots\wedge\mu(\Sigma_l)\right).
\]
Hence
\begin{eqnarray*}
\int_{\mathcal{M}}\mu(\Sigma_1)\wedge\ldots\wedge\mu(\Sigma_l)
&=&\lim_{\eps\mapping 0}\int_{\mathcal{M}_\eps}\mu(\Sigma_1)\wedge\ldots\wedge\mu(\Sigma_l)\\
&=&\lim_{\eps\mapping 0}\int_{\partial\mathcal{M}_\eps}\mu(\Xi_1)\wedge\mu(\Sigma_2)\wedge\ldots\wedge\mu(\Sigma_l)\\
&=&\int_{\sphere{4}}\tau^*\mu(\Xi_1)\wedge\tau^*\mu(\Sigma_2)\wedge\ldots
\wedge\tau^*\mu(\Sigma_l)\\
& &\hbox{for some appropriate tractator $\tau$.}\\
&=&\int_{\sphere{4}}\mathrm{PD}(\Xi_1)\wedge\mathrm{PD}(\Sigma_2)\wedge\ldots
\wedge\mathrm{PD}(\Sigma_l)\\
&=&\hbox{intersection number of }\bigcap_{i=2}^{l}\Sigma_i\hbox{ with }\Xi_1\\
&=&\hbox{linking number of }\Sigma_1 \hbox{ and }\bigcap_{i=2}^{l}\Sigma_i
\end{eqnarray*}
The possible configurations for Donaldson numbers in the case $k=1$ depend on
the relation
\[
\sum_{j=1}^{l}d_i=4l-8k+3=4l-5
\]
where $d_i$ is the dimension of the submanifold $\Sigma_i$. When $l=2$, we know
that $d_1+d_2=3$, so the only possible configurations are
\[
\begin{array}{c|c}
 d_1 & d_2 \\ \hline\hline
 0   & 3\\ \hline
 1& 2 \\ \hline
\end{array}.
\]
\begin{center}
$\epsfig{file=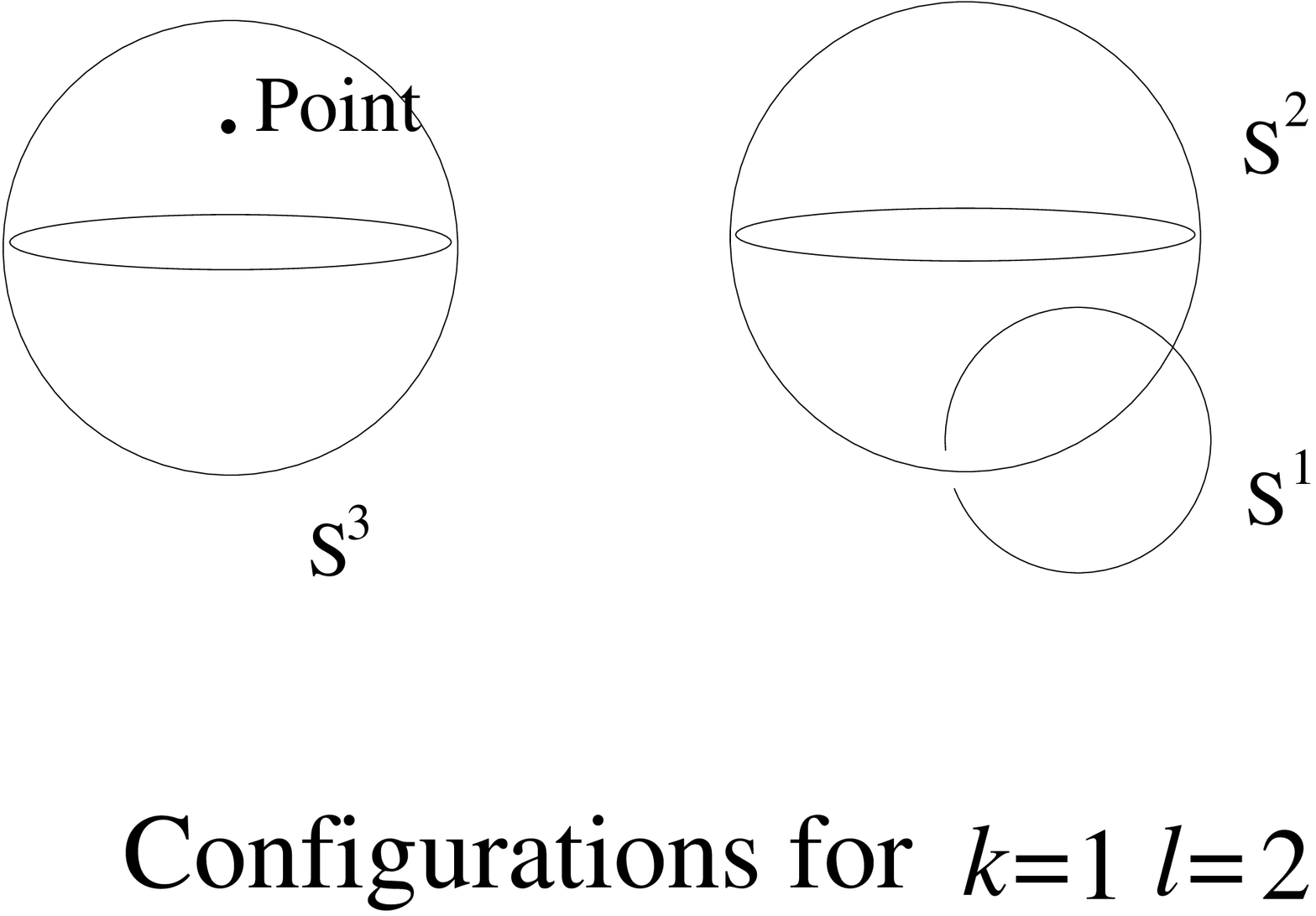,width=7cm,angle=-90}$
\end{center}
For $l=3$ we have the following linking configurations.
%\begin{figure}[!ht]
%\input{Configk1.pstex_t}
\begin{center}
$\epsfig{file=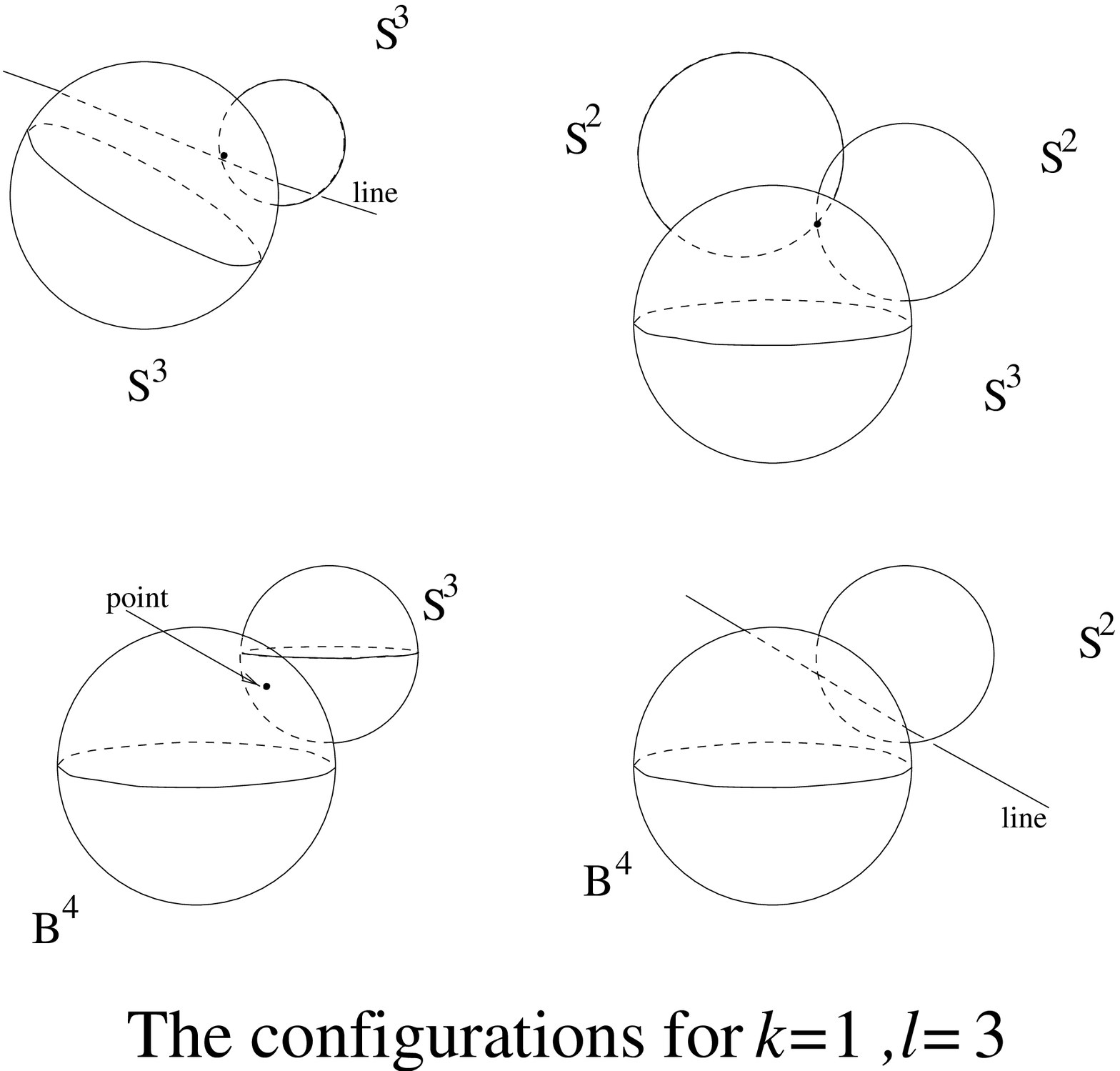,width=12cm}$
\end{center}
Our next task will be to examine $k\ge2$.

%\subsection{Higher $k$}

\subsection{Is there any linking for $k\ge 2$?}
It would be prudent to examine $\P(\Sigma_1,\ldots,\Sigma_l)$.
Recall that
\begin{eqnarray*}
& &\P(\Sigma_1,\ldots,\Sigma_l)(\phi)\\
&=&\int_{\quat^{k-1}\times{\mathcal{M}}_1}
\iota^*\phi\,\left(\mu^1(\Sigma_1)+\sum_{j=1}^{k-1}p_j^*PD(\Sigma_1)\right)\wedge
\ldots\wedge\left(\mu^1(\Sigma_l)+\sum_{j=1}^{k-1}p_j^*PD(\Sigma_l)\right)
\end{eqnarray*}
and let us look at the various configurations. Choose pairwise disjoint,
submanifolds
$\Sigma_1,\ldots\Sigma_l$ of $\Rf$ of dimensions $d_1,\ldots,d_l$
respectively. Then
\begin{equation}
\sum_{i=1}^l(4-\d_i)=8k-3
\end{equation}
yields
\begin{equation}
\sum_{i=1}^l d_i=4l-8k+3.
\end{equation}
\begin{defn1}
Let $K$ be a finite set. Then
\begin{enumerate}
\item we call an $n$-tuple $I$ of elements of $K$ an ordered subset of $K$
whenever $I=(i_1,\ldots,i_p)$ we have $i_j\neq i_r$ for all $j\neq
r$;
\item we say $n\in I=(i_1,\ldots,i_p)$ if there is $j$ such that $i_j=n$;
\item we shall write $\#I=p$ if $I=(i_1,\ldots i_p)$;
\item we shall say that a collection $I_1,\ldots I_m$ of ordered subsets of
$K$ is a partition of $K$ if for each $n\in K$ there is precisely
one $r$ such that $n\in I_r$.
\end{enumerate}
\end{defn1}
Now as we said earlier, each term in $\P$ is
$\pm\mathrm{Don}_1(\Sigma_{I_1})\prod_{j>1}^p\sharp(\Sigma_{I_i})$
for ordered subsets
 $I_1,\ldots,I_p$ that partition $\{1,\ldots,l\}$. For this to
give a nonzero contribution to $\P$ we need $I_1$ to have at least
2 elements and each of the other $I_i$ at least 1 element. Now
$I_1$ must satisfy
\begin{equation}
\sum_{j\in I_1}(4-\d_j)=8(1)-3=5
\end{equation}
that is
\begin{equation}
\sum_{j\in I_1}d_{j}=4\#I_1-5.
\end{equation}
%and so $p\ge 2$.
Now we know that for each of $I_j,j>1$ we need
$\sharp(\Sigma_{I_j})\neq 0$ so we must have for each $j>2$,
\[
\sum_{q\in I_j}(4-d_{q})=4
\]
that is
\[
\sum_{q\in I_j}(d_{q})=4\#I_j-4.
\]
Now
\begin{eqnarray*}
\sum_{j=1}^{l}d_j&=&\sum_{n=1}^p\sum_{j\in I_n} d_j\\
&=&\sum_{j\in I_1} d_j +\sum_{n=2}^p\sum_{j\in I_n} d_j\\
&=&4\#I_1-5+\sum_{n=2}^p(4\#I_n-4)\\
&=&4\sum_{n=1}^p\#I_n-5-4(p-1)\\
&=&4l-4p-1.
\end{eqnarray*}
Thus we must have
\begin{eqnarray*}
4l-4p-1=4l-8k+3
\end{eqnarray*}
i.e.
\begin{eqnarray*}
p=2k-1.
\end{eqnarray*}
So $k$ controls the number of ordered subsets of $\{1,\ldots,l\}$
that form a partition, moreover this number has to be odd.\par
However, if we examine the form of $\P$ more closely,
\begin{eqnarray*}
& &\P(\Sigma_1,\ldots,\Sigma_l)(\phi)\\
&=&\int_{\quat^{k-1}\times{\mathcal{M}}_1}\iota^*\phi
\left(\mu^1(\Sigma_1)+\sum_{j=1}^{k-1}p_j^*PD(\Sigma_1)\right)\wedge
\ldots\wedge\left(\mu^1(\Sigma_l)+\sum_{j=1}^{k-1}p_j^*PD(\Sigma_j)\right),
\end{eqnarray*}
we see that each term is the product of sums of $k$ terms, so for
 any ordered subsets $I_1,\ldots,I_p$  that partition $\{1,\ldots,l\}$ and give non-zero contribution to $\P$ must satisfy
\[
p\leq k.
\]
Thus we have
\[
k\ge p=2k-1
\]
which is impossible for $k>1$. Hence $\P$ is a trivial
topological distribution, assigning 0 to any set of pairwise disjoint, compact
submanifolds of
$\Rf$ and any test function $\phi$. We have therefore proved
\begin{thm1}\label{conclusion}
For $k\ge2$, and any compactly supported $\phi$ there are no
anomalies, i.e for $k\ge2$
\[
\mathrm{Don}_k(\Sigma_1,\ldots,\Sigma_l)(\phi)=0
\]
for all pairwise disjoint, compact submanifolds $\Sigma_1,\ldots,\Sigma_l$.
\end{thm1}
We proved this for the resolution of the moduli space. Since the
integral is identically 0 on the resolution, it must also be zero
on the moduli space itself, thus agreeing with the infinite
dimensional construction.
%\begin{exm1} We examine the case $k=2,l=4$. So the number of possible ordered
%subsets that partition $\{1,\ldots, 4\}$ and make a non-zero
%contribution to $\P$ is $3$, which we shall call $I_1,I_2,I_3$.
%For
%\[
%\mathrm{Don_1}(\Sigma_{I_1})\neq 0
%\]
%we need $\#I_1\ge 2$. We also need $\#I_2,\#I_3\ge 1$, and
%$\#I_1+\#I_2+\#I_3=4$, and together this implies that $\#I_1=2$
%and $\#I_2=\#I_3=1$. Now recall that for $j=2,3$
%\[
%\sum_{q\in I_j}(d_{q})=4\#I_j-4.
%\]
%This implies that $\Sigma_{I_j}$ is a finite collection of $r_j$
%points for $j=2,3$. Thus if $I_1=\{m_1,m_2\},
%I_2=\{m_3\},I_4=\{m_4\}$ the contribution that this makes to $\P$
%is
%\[
%\mathrm{Don}_1(\Sigma_{m_1},\Sigma_{m_2})\sharp(\Sigma_{m_3})\sharp(\Sigma_{m_3})=r_2r_3\mathrm{link}(\Sigma_{m_1},\Sigma_{m_2}).
%\]
%We have already said that for $k=1,l=2$ there are only two
%possibilities for the dimensions of $\Sigma_{m_1}$ and
%$\Sigma_m_2$, namely that
%\end{exm1}
\subsection{Concluding Remarks}
Although we have had something of a disappointment that there is no linking
number for $k>1$ on the moduli space of instantons, nor on any resolutions, we
have developed some potentially
powerful techniques in finding formul{\ae} for the cohomology of a
hyperK\"ahler reduction. One can
hope that the technique for hyperK\"ahler manifolds with boundaries may produce
information about the topology of the higher instanton spaces by looking at the
topology of the end and concluding that the Moduli space is a cone on this
manifold. Also there may be something to be said about the perturbed moduli
spaces with their relationship with the Seiberg-Witten equations.
% The
%techniques introduced here may also provide valuable information about the
%non-Abelian Seiberg-Witten equations if we can compactify the gauge group.
%\par
%We have also found that on nearby resolutions of the moduli space, there
%is a very complicated linking theory involving intersection numbers and linking
%numbers.
\par
The Author would like to thank the EPSRC for their financial support, and
is appreciative of (in no particular order) for the very welcome input of Prof Dr. Dr. Victor Pidstrigatch, Dr. Mario Micallef,  Prof. John Rawnsley, Dr. Roger Bielawski and Dr. Richard Thomas. HOAMGD.
%\backmatter
\bibliography{Bibtex}

\begin{thebibliography}{10}

\bibitem{An}
{\sc Anselmi, D.}
\newblock Anomalies in instanton calculus.
\newblock {\em Nuclear Phys. B 439}, 3 (1995), 617--649.

\bibitem{An2}
{\sc Anselmi, D.}
\newblock Topological field theory and physics.
\newblock {\em Classical Quantum Gravity 14}, 1 (1997), 1--20.

\bibitem{ADHM}
{\sc Atiyah, M.~F., Hitchin, N.~J., Drinfel'd, V.~G., and Manin, Y.~I.}
\newblock Construction of instantons.
\newblock {\em Phys. Lett. A 65}, 3 (1978), 185--187.

\bibitem{BGV}
{\sc Berline, N., Getzler, E., and Vergne, M.}
\newblock {\em Heat kernels and {D}irac operators}.
\newblock Springer-Verlag, Berlin, 1992.

\bibitem{DO}
{\sc Donaldson, S.~K.}
\newblock Connections, cohomology and the intersection forms of $4$-manifolds.
\newblock {\em J. Differential Geom. 24}, 3 (1986), 275--341.

\bibitem{DK}
{\sc Donaldson, S.~K., and Kronheimer, P.~B.}
\newblock {\em The geometry of four-manifolds}.
\newblock The Clarendon Press, Oxford University Press, New York, 1990.
\newblock Oxford Science Publications.

\bibitem{MU2}
{\sc Munn, J.}
\newblock Equivariant {I}ntegration {F}ormul{\ae} in {H}yper{K}\"ahler
  {G}eometry.
\newblock {\em (To Appear)\/}.

\bibitem{MU}
{\sc Munn, J.}
\newblock The {ADHM} construction and its applications to {D}onaldson {T}heory.
\newblock {\em Thesis for the Degree of Ph.D at the University of Warwick\/}
  (2001).

\bibitem{N}
{\sc Nahm, W.}
\newblock Self-dual monopoles and calorons.
\newblock In {\em Group theoretical methods in physics (Trieste, 1983)}.
  Springer, Berlin, 1984, pp.~189--200.

\bibitem{S}
{\sc Selby, M.}
\newblock Donaldson invariants and equivariant cohomology.
\newblock {\em Thesis for the degree of D.Phil at the University of Oxford\/}
  (1998).

\bibitem{WO}
{\sc Wood, R. M.~W.}
\newblock Quaternionic eigenvalues.
\newblock {\em Bull. London Math. Soc. 17}, 2 (1985), 137--138.

\end{thebibliography}
 \bibliographystyle{acm}
\end{document}